\pdfoutput=1
\documentclass[a4paper,fleqn]{cas-sc-upload}

\usepackage[numbers,sort&compress]{natbib}

\usepackage{placeins}
\usepackage{tabularx}
\usepackage{needspace}
\usepackage[ruled,lined]{algorithm2e}
\usepackage{doi} % Make every bibliography DOI a clickable https://doi.org link.

\renewcommand{\sssectionfont}{\rmfamily\fontsize{10.5pt}{12pt}\bfseries\upshape\selectfont}
\numberwithin{equation}{section}
\graphicspath{{./}}

\makeatletter
\let\cas@emails\printemails
\let\cas@orcid\printorcid
\let\cas@fnotes\printfnotes
\let\printemails\relax
\let\printorcid\relax
\def\printfnotes{\cas@fnotes\cas@orcid\cas@emails}
\makeatother

\ExplSyntaxOn
\cs_set:Npn \__first_footerline:
  {
    \group_begin:
    \small
    \sffamily
    \ifnum\theblind>0\relax
    \else
      \__short_authors: :~
    \fi
     \rmfamily \itshape Preprint 
    \group_end:
  }
\ExplSyntaxOff

\begin{document}
\let\WriteBookmarks\relax
\def\floatpagepagefraction{1}
\def\textpagefraction{.001}

\shorttitle{Error attribution in GRW methods}
\shortauthors{S. Abkin and P. Daripa}

\title[mode=title]{\texorpdfstring{Error attribution in gradient random walk methods for parabolic equations}{Error attribution in gradient random walk methods for parabolic equations}}

\author[1]{Stephen Abkin}[orcid=0009-0000-5981-9721]
\cormark[1]
\ead{stephen122204@tamu.edu}
\credit{Conceptualization, Methodology, Software, Validation, Formal analysis, Investigation, Data curation, Visualization, Writing -- original draft, Writing -- review and editing}

\author[1]{Prabir Daripa}[orcid=0000-0002-8771-0149]
\ead{daripa@tamu.edu}
\credit{Conceptualization, Methodology, Supervision, Project administration, Writing -- review and editing}

\affiliation[1]{organization={Department of Mathematics, Texas A\&M University},
                city={College Station},
                postcode={77843},
                state={TX},
                country={USA}}

\cortext[1]{Author to whom correspondence should be addressed.}

\begin{abstract}[Abstract]

Gradient random walk (GRW) methods represent a solution's spatial derivative by weighted particles and recover the solution by cumulative summation. We examine their accuracy for the one-dimensional heat, scalar cubic reaction--diffusion, and Burgers' equations. Thirty-seed ensembles, paired reconstructions, and deterministic controls identify the sources of two apparent accuracy limits. For the heat equation, comparing cumulative bin sums at bin centers introduces a first-order alignment error, which the exact reflected-particle distribution predicts and bin-edge comparison removes. An identity for the expected squared error predicts the remaining floor and the particle count for a target root-mean-square accuracy. With unequal weights, a count chosen in advance gives an error within $1\%$ of its prediction. Reaction-front position, speed, and shape errors decrease under refinement. For Burgers' equation, exact transformed boundary data reduce the deterministic Cole--Hopf error about $490$-fold, while fixed-endpoint particle runs approach their deterministic reference.

\end{abstract}

\begin{keywords}[MSC 2010]
\raggedright Primary 65C35\sep Secondary 65C05, 65M75, 65M15, 35K57\par
\par\medskip\noindent\textit{Keywords}:\par
gradient random walk \sep stochastic particle methods \sep numerical error estimation \sep cumulative reconstruction \sep reaction--diffusion equations \sep Cole--Hopf transformation
\end{keywords}

\maketitle

\section{Introduction}
\label{sec:intro}

Gradient random walk (GRW) methods approximate parabolic equations by transporting their spatial derivative on weighted particles. Independent Gaussian displacements represent diffusion, and cumulative summation of the particle weights recovers the solution~\citep{chorin1973,ghoniem1985}. The measured error also depends on where that reconstruction is compared with a reference and how the boundary data are imposed. Identifying these effects is essential when choosing a particle count or interpreting an error plateau.

We separate these sources of error for the heat equation, a scalar reaction--diffusion equation with a FitzHugh--Nagumo-type cubic reaction term, and Burgers' equation. We use standard test problems with known exact solutions to assess the accuracy of the calculations. Their one-dimensional setting allows us to isolate individual numerical operations using analytical predictions, exact particle distributions for the heat equation, and deterministic reference calculations. We run thirty-seed ensembles to quantify sampling variability and compare paired reconstructions while holding the particle trajectories fixed. To evaluate the effects of reconstruction and boundary choices, we perform deterministic controls using a known particle distribution or a finite-difference solution. We interpret these comparisons through the principles of code verification~\citep{roache1998,oberkampfroy2010}.

The main results are as follows:
\begin{itemize}
\item For the heat equation, a half-bin displacement between a cumulative reconstruction and its comparison points explains the observed error floor at fixed bin count. Comparing at the represented bin edges removes this alignment error. An exact identity for the expected squared error then predicts the remaining floor, the transition from sampling error to that floor, and a particle count for a specified root-mean-square accuracy. A second heat problem with unequal particle weights tests these predictions on configurations selected before the simulations, including an error prediction confirmed to within $1\%$ on fresh seeds in Section~\ref{sec:heat-two-step}.
\item For the reaction--diffusion equation, we follow the same diffusion step with a deterministic update of the weights. We verify the combined calculation against an exact traveling front and measure decreasing errors in its profile, location, speed, and aligned shape in Section~\ref{sec:fhn-convergence}.
\item For Burgers' equation, we use controlled comparisons to identify the transformed boundary data as the source of the dominant deterministic error. We reduce this error by a factor of about $490$ by supplying exact data to the deterministic Cole--Hopf solve. With the original fixed endpoints, the particle calculation approaches the corresponding deterministic reference under initialization-point refinement in Section~\ref{sec:colehopf-decoupled}.
\end{itemize}

These results build on earlier analyses of particle accuracy. Ghoniem and Sherman~\citep{ghoniem1985} developed finite-interval reflection rules and compared the reconstructed fields with analytical solutions. Mascagni~\citep{mascagni1995} separated wave-speed and spatial-discretization errors for a deterministic Nagumo particle analogue. Bertaglia, Pareschi, and Caflisch~\citep{bertaglia2024} separated sampling and mesh errors for gradient-based Monte Carlo methods in a hyperbolic relaxation setting. They give the exact pointwise root-mean-square sampling error, and their gradient-particle formulation includes signed particles of equal mass magnitude. Here, we combine the classical sampling variance with the reflected finite-interval distribution, the comparison-location error in the reported discrete norm, and unequal weight magnitudes with fixed allocation. We use the resulting predictions to guide the choice of comparison points and particle count. For Burgers' equation, we examine the additional boundary and recovery operations introduced by the Cole--Hopf transformation.

The three equations introduce these operations successively:
\[
\begin{aligned}
    u_t &= \alpha u_{xx},\\
    u_t &= D u_{xx}+f(u),\\
    u_t+u u_x &= \nu u_{xx}.
\end{aligned}
\]
Here \(\alpha\) is the thermal diffusivity, \(D\) is the diffusion coefficient, \(f(u)\) is the cubic reaction term specified in Section~\ref{sec:fhn-benchmark}, and \(\nu\) is the kinematic viscosity. We first test diffusion and cumulative reconstruction with the heat equation, then include the weight update in the scalar reaction--diffusion calculation. For Burgers' equation, we differentiate the transformed heat solution and apply the nonlinear recovery. Section~\ref{sec:benchmarks} specifies the initial and boundary data and the exact reference solutions.

The underlying particle construction has a classical history. Chorin~\citep{chorin1973} represented viscous diffusion through Gaussian random displacement in the random-vortex method. Ghoniem and Sherman~\citep{ghoniem1985} explained how integration reduces fluctuations in the recovered field and gave a deterministic reaction-weight update. Sherman and Peskin~\citep{shermanpeskin1986} studied a related reaction--diffusion construction with random creation and removal of particles, and later applied related ideas to the Hodgkin--Huxley cable equations~\citep{shermanpeskin1988}. We use the deterministic multiplier of Ghoniem and Sherman, evaluated after diffusion as specified in Section~\ref{sec:reaction-update}. Puckett~\citep{puckett1989} proved convergence for a related random particle method for the Kolmogorov--Petrovskii--Piskunov equation, with bound \(O((\ln N)N^{-1/4})\) under \(\Delta t=O(N^{-1/4})\), and observed decay approximately proportional to \(N^{-1/2}\). Sherman and Mascagni~\citep{shermanmascagni1994} extended gradient random walk to two dimensions and demonstrated its concentration of particles where the field gradient is large.

For Burgers' equation, the transformations of Hopf~\citep{hopf1950} and Cole~\citep{cole1951} reduce the problem to the heat equation. Roberts~\citep{roberts1989} proved convergence of a direct fractional-step random-walk method that handles the nonlinear term by particle advection. Aref and Daripa~\citep{arefdaripa1984} used Burgers' equation to examine spurious behavior caused by failure to conserve discrete kinetic energy in finite-difference schemes. Yan~\citep{yan2023} treated nonhomogeneous Burgers problems on bounded domains through a generalized Cole--Hopf transformation with derivative boundary conditions. These studies motivate checking both the transformed boundary data and the recovery of the physical field in the particle calculation considered here.

We use \emph{globs} for the weighted gradient particles. Their cumulative reconstruction is defined at every position, with jumps at the particle locations~\citep{puckett1989,bertaglia2024}. We report it at specified \emph{reconstruction points}, and call intervals used to group weights \emph{reconstruction bins}. These points and bins serve field evaluation and error measurement, while the particles move continuously in space. The abbreviation GRW denotes gradient random walk throughout. The lattice-based generalized or global random walk~\citep{vamos2003,suciu2021} and the pointwise elliptic walk-on-spheres and walk-on-stars solvers~\citep{sawhney2020,sawhneymiller2023} are distinct methods.

Section~\ref{sec:method} presents the particle formulation, and Section~\ref{sec:benchmarks} specifies the three canonical problems. Section~\ref{sec:verification} defines the reference comparisons, error measures, and ensemble analysis. Sections~\ref{sec:heat}--\ref{sec:colehopf} report the numerical studies. Section~\ref{sec:analysis} derives the identities and compatibility conditions that determine their applicability to other problems. Section~\ref{sec:results-summary} collects the quantitative findings, and Section~\ref{sec:discussion} concludes.
\section{Gradient Random Walk Formulation}
\label{sec:method}
 
\subsection[Particle representation of the gradient]{Particle representation of the gradient}
 
Let \(v=u_x\) denote the spatial derivative of the solution field. In one dimension, GRW represents~\(v\) by a collection of \(N\)~weighted point masses:
\begin{equation}
    v_N(x,t)=\sum_{i=1}^{N} w_i(t)\,
    \delta\bigl(x-X_i(t)\bigr),
    \label{eq:gradient-measure}
\end{equation}
where \(X_i(t)\) is the position of the \(i\)-th glob, \(w_i(t)\) is its signed weight, and \(\delta\) is the Dirac delta function. Each glob weight contributes to the spatial derivative~\(u_x\). This gradient representation and cumulative recovery follow Ghoniem and Sherman~\citep{ghoniem1985} and are summarized in the modern gradient-based Monte Carlo framework of Bertaglia et al.~\citep{bertaglia2024}.
 
To see how the initialization works, consider a step initial condition with a jump from~\(u_L\) to~\(u_R\) at position~\(x_0\). At a jump discontinuity, the derivative~\(u_x\) is concentrated entirely at the jump location and has total value~\(u_R - u_L\). We represent this by placing all~\(N\) globs at the jump location with equal weights given by
\begin{equation}
    X_i^0=x_0,\qquad
    w_i^0=\frac{u_R-u_L}{N},
    \qquad i=1,\ldots,N,
    \label{eq:jump-initialization}
\end{equation}
so that the weights sum to~\(\sum_i w_i^0=u_R-u_L\), matching the total jump. For smooth monotone initial profiles such as the logistic front used in the reaction--diffusion problem, the globs are placed at positions that spread the total weight evenly across the profile rather than concentrating it at a single point. We describe this construction in more detail in Section~\ref{sec:fhn-benchmark}.
 
\subsection{Brownian displacement and diffusion}

\label{sec:brownian}

At each time step, every glob is independently displaced by a Gaussian random increment. For a diffusion coefficient \(\kappa>0\), the update rule is~\citep{chorin1973,ghoniem1985}
\begin{equation}
    \widetilde X_i^{\,n+1}
    =X_i^n+\sqrt{2\kappa\Delta t}\,Z_i^n,
    \qquad
    Z_i^n\sim\mathcal N(0,1),
    \label{eq:brownian-update}
\end{equation}
where \(X_i^n\) is the position of glob~\(i\) at time~\(t_n=n\Delta t\), \(\Delta t\) is the time step, and \(Z_i^n\) is a standard normal random variable drawn independently for each glob and each step. The quantity~\(\widetilde X_i^{\,n+1}\) is the glob position after the Brownian displacement and before boundary reflection. The diffusion coefficient takes the value \(\kappa=\alpha\) for the heat equation, \(\kappa=D\) for the scalar reaction--diffusion equation, and \(\kappa=\nu\) for the transformed Burgers heat equation.
 
In the absence of boundaries, the variance of the accumulated displacement over time~\(t\) is~\(2\kappa t\), which is the variance of the fundamental solution of~\(u_t = \kappa u_{xx}\)~\citep{chorin1973,ghoniem1985}.
 
\subsection{Cumulative reconstruction}
\label{sec:reconstruction}
 
Once the globs have been displaced, we recover the solution field~\(u\) by sorting the globs in order of increasing position and summing their weights up to each evaluation point. Given a known left-boundary value~\(u_L\), the reconstructed field at any point~\(x\) is
\begin{equation}
    u_N(x)=u_L+\sum_{i:X_i\leq x}w_i.
    \label{eq:cumulative-reconstruction}
\end{equation}
Sorting lets us evaluate these sums efficiently by following the globs in spatial order. Each sum includes exactly the globs with \(X_i\leq x\). In practice, the partial sums are evaluated at the reconstruction points introduced in Section~\ref{sec:intro}. We take these to be the finite uniform set~\(\{x_j\}_{j=0}^{M-1}\) with spacing~\(h\), independent of the glob positions. At each parameter value, all realizations use the same reconstruction points. A refinement study may also vary these points, as specified in its design. Throughout this paper, \(M\)~denotes the number of reconstruction locations, and \(N\)~the number of globs for the heat and reaction--diffusion problems. Two reconstruction operations appear. In the first, the cumulative sum is evaluated directly at \(M\) uniform points, as for the scalar reaction--diffusion problem and the coupled heat treatment of Section~\ref{sec:heat-ensemble}. In the second, the glob weights are first assigned to \(M\) bins and then summed in order across the bins, as for the fixed-bin heat treatments and the Cole--Hopf Burgers construction of Section~\ref{sec:burgers-benchmark}. All comparisons with exact or reference solutions use the reconstructed field at these locations.

\subsection{Grid-free evolution and the role of grids}
\label{sec:grid-free}

The globs always move through continuous space by the Brownian displacement~\eqref{eq:brownian-update}, never on a lattice. For the heat and scalar reaction--diffusion problems the entire time evolution is grid free. In particular, the reaction-weight update reads the field at the glob positions themselves using the cumulative reconstruction~\eqref{eq:cumulative-reconstruction}, so no spatial grid is consulted during a run. The reconstruction points enter only after the evolution, as the measurement device that Section~\ref{sec:norms} uses to define the error.

The Cole--Hopf Burgers calculation is the exception. Because the random walk is applied to the transformed heat variable~\(\phi\), the construction itself discretizes space at fixed locations. \(P\) initialization points build the transformed initial gradient and place the \((P-1)\) \(\phi_x\)-globs, and the transported weights are then assigned to \(M\) reconstruction bins, smoothed, and cumulatively summed to reconstruct \(\phi\). We differentiate this reconstruction and recover \(u=-2\nu\phi_x/\phi\). These discretizations belong to the Cole--Hopf construction rather than to the gradient random walk, and Section~\ref{sec:colehopf-decoupled} varies~\(P\) and~\(M\) independently. The Brownian displacement of the~\(\phi_x\)-globs remains grid free.

\subsection{Boundary reflection}
\label{sec:boundary-reflection}

On a finite computational interval~\([0,L]\), any glob that crosses the left boundary at~\(x=0\) is reflected back by the map \(X_i\mapsto-X_i\), and any glob that crosses the right boundary at~\(x=L\) is reflected by \(X_i\mapsto2L-X_i\). In both cases, the reflection preserves the distance by which the glob overshot the boundary rather than placing the glob directly on the wall. After any required boundary reflection, we denote the resulting position by \(X_i^{*}\in[0,L]\). No later operation in the step moves the globs, so \(X_i^{n+1}=X_i^{*}\). The reflection is applied repeatedly until the position lies in $[0,L]$, because a Gaussian increment is unbounded. For constant diffusivity this folded Gaussian step is the reflected Brownian transition kernel on $[0,L]$, and its composition over the steps of a run is the reflected kernel at every step time, so the time step partitions the exact reflected evolution without discretizing it. The distribution of a particle position is then the classical distribution of reflected Brownian motion on $[0,L]$, whose density is given by the method of images~\citep{carslawjaeger1959,pavliotis2014} and is written out in Appendix~\ref{app:heat-uncertainty}. We denote its distribution function by $F_R$ and use it as an exact reference in Section~\ref{sec:heat-grid-paired}. Section~\ref{sec:analysis-cdf} gives the argument that the folded steps sample it exactly.
 
Finite-domain boundary treatment for gradient particles follows the image and reflection construction of Ghoniem and Sherman~\citep{ghoniem1985}. Under the particle representation of the gradient, the way the signed weight is handled upon reflection depends on the type of boundary condition. For Dirichlet boundaries (i.e., prescribed boundary values of~\(u\)), the weight is preserved upon reflection. For a reaction--diffusion equation this represents a fixed boundary state only when that state is an equilibrium of the reaction, as shown in Section~\ref{sec:analysis-reaction}. For homogeneous Neumann boundaries (i.e., zero-flux conditions on~\(u\)), the weight is negated upon reflection. A preserved weight extends~\(u_x\) evenly across the wall, and a negated weight extends it oddly, which enforces~\(u_x=0\) at a homogeneous Neumann boundary.
 
In the Cole--Hopf Burgers calculation, GRW evolves globs representing~\(\phi_x\), where~\(\phi\) satisfies the heat equation~\(\phi_t=\nu\phi_{xx}\) on~\([0,L]\). We preserve the glob weights upon reflection, reconstruct~\(\phi\) using~\(\phi_0(0)\) as the integration constant, and correct the weight sum to match the prescribed endpoint difference~\(\phi_0(L)-\phi_0(0)\). These operations specify the finite-domain boundary treatment for~\(\phi\). The boundary values of~\(u=-2\nu\phi_x/\phi\) then follow from the reconstructed transformed field.
 
\subsection{Reaction-weight update}
\label{sec:reaction-update}
 
For reaction--diffusion equations of the form
\[
    u_t=D u_{xx}+f(u),
\]
differentiating with respect to~\(x\) yields an evolution equation for the derivative~\(v = u_x\):
\[
    v_t=Dv_{xx}+f'(u)v.
\]
The first term on the right,~\(Dv_{xx}\), is a diffusion term and is handled by the Brownian displacement step~\eqref{eq:brownian-update}. The second term,~\(f'(u)v\), describes how the reaction modifies the derivative, and is handled by updating the particle weights after the diffusion step. The diffusion and reflection steps produce intermediate positions~\(X_i^{*}\) and weights~\(w_i^{*}\), from which we reconstruct the intermediate field
\[
    u_i^{*}=u_L+\sum_{j:\,X_j^{*}\leq X_i^{*}} w_j^{*},
\]
and then apply the deterministic weight update of Ghoniem and Sherman~\citep{ghoniem1985}:
\begin{equation}
    w_i^{n+1}=w_i^{*}\left[1+\Delta t\,R(u_i^{*})\right],
    \qquad R(u)=f'(u).
    \label{eq:reaction-weight-update}
\end{equation}
The star specifies evaluation after diffusion and reflection but before the reaction update, so each step proceeds as \((X_i^n,w_i^n)\to(X_i^{*},w_i^{*})\to u_i^{*}\to w_i^{n+1}\). Sherman and Peskin~\citep{shermanpeskin1986} treat the reaction step as a random birth--death event, duplicating a glob with probability \(\Delta t\,R(u)\) when \(R(u)>0\) and removing it with probability \(-\Delta t\,R(u)\) when \(R(u)<0\), provided $\Delta t\,|R(u)|\le1$. The deterministic multiplier \(1+\Delta t\,R(u_i^{*})\) in~\eqref{eq:reaction-weight-update} is the conditional mean of that event given the current state. This conditional-mean identity applies to one step at a given reconstructed state. The evolving field couples the globs, so we assess the full deterministic calculation through the ensembles of Section~\ref{sec:fhn-convergence}. Section~\ref{sec:analysis-reaction} discusses in more detail the class of reaction laws and boundary data to which the update applies.
 
Diffusion moves the globs, while the reaction update changes their contributions to the gradient. The size of the multiplier for the reaction--diffusion front is quantified in Section~\ref{sec:fhn-representative}.
 
Algorithm~\ref{alg:grw} summarizes the GRW time-stepping procedure. For the heat equation, there is no reaction term in the PDE, so the weight update is omitted. For the scalar reaction--diffusion equation, the weight update is included, with \(R(u)\) given by Eq~\eqref{eq:fhn-statistic} of Section~\ref{sec:fhn-benchmark}. Changing the reaction also requires compatible initial data, boundary treatment, and time-step selection, as discussed in Section~\ref{sec:analysis-reaction}.
 
\FloatBarrier
\begin{algorithm}[!htbp]
\caption{GRW time-stepping for the heat and scalar reaction--diffusion problems (reaction term optional)}
\label{alg:grw}
\KwIn{diffusivity \(\kappa\), domain \([0,L]\) with boundary type (Dirichlet or Neumann), time step \(\Delta t\), left state \(u_L\), number of globs \(N\), number of steps \(K\), weight update function \(R(u)\) (optional).}
\vspace{0.3em}
\KwOut{reconstructed field \(u_N\) at the reconstruction points at time \(T = K\Delta t\), with the sorted glob positions and weights.}
\vspace{0.5em}
Initialize positions \(X_i^0\) and weights \(w_i^0\) via Eq~\eqref{eq:jump-initialization} or Eq~\eqref{eq:fhn-quantile-initialization} of Section~\ref{sec:fhn-benchmark}. \\
\vspace{0.3em}
\For{\(n = 0, \ldots, K-1\)}{
    \textit{Diffuse:} displace each glob via Eq~\eqref{eq:brownian-update}. \\
    \textit{Reflect:} apply the boundary maps of Section~\ref{sec:boundary-reflection}, negating the weight at Neumann boundaries, giving intermediate positions \(X_i^{*}\) and weights \(w_i^{*}\). \\
    \textit{Sort} globs by position, the weights following the same ordering. \\
    \textit{Reconstruct} \(u_i^{*}\) by cumulatively summing the intermediate weights in position order. \\
    \If{reaction term present}{
        Update weights via Eq~\eqref{eq:reaction-weight-update}, with $R(u)=f'(u)$. \\
    }
}
\vspace{0.25em}
Evaluate \(u_N\) at the reconstruction points via Eq~\eqref{eq:cumulative-reconstruction}. \\
\vspace{0.15em}
\end{algorithm}
\FloatBarrier

\section{Canonical Problems}
\label{sec:benchmarks}
 
\subsection{Heat equation step problem}
\label{sec:heat-benchmark}
 
In the heat problem, we solve the diffusion equation
\begin{equation}
    u_t=\alpha u_{xx},
    \qquad x\in[0,L],
    \label{eq:heat-pde}
\end{equation}
with a step initial condition that jumps from~\(u_L\) to~\(u_R\) at position~\(x_0\),
\begin{equation}
    u(x,0)=
    \begin{cases}
        u_L, & x<x_0,\\
        u_R, & x>x_0.
    \end{cases}
    \label{eq:heat-step-data}
\end{equation}
For the corresponding problem on the real line, convolution of the step data with the heat kernel gives the classical error-function profile~\citep{carslawjaeger1959}
\begin{equation}
    u^{\mathrm{ex}}(x,t)
    =u_L+\frac{u_R-u_L}{2}
    \left[
      1+\operatorname{erf}
      \left(\frac{x-x_0}{2\sqrt{\alpha t}}\right)
    \right].
    \label{eq:heat-exact}
\end{equation}
The globs are initialized at the jump location via~\eqref{eq:jump-initialization}, displaced by~\eqref{eq:brownian-update} with \(\kappa=\alpha\), and the field is recovered by~\eqref{eq:cumulative-reconstruction}. The calculation therefore tests Brownian displacement and cumulative reconstruction without the reaction update or inverse transformation used in the later problems.
 
\subsection{Scalar reaction--diffusion traveling front}
\label{sec:fhn-benchmark}
 
The second canonical problem is a scalar reaction--diffusion equation
\begin{equation}
    u_t=D u_{xx}+f(u),
    \label{eq:fhn-pde}
\end{equation}
related to the excitable-nerve models introduced by FitzHugh~\citep{fitzhugh1961} and Nagumo et al.~\citep{nagumo1962}.

We use a FitzHugh--Nagumo-type cubic reaction term, with the coefficients and parameter regime specified below. We verify the reaction-driven weight update combined with diffusion against the exact traveling front

\begin{equation}
    u^{\mathrm{ex}}(x,t)
    =\frac{1}
    {1+\exp\left[-(x+\theta t-x_c)/2\right]},
    \qquad
    \theta=\sqrt{2}\,(0.5-a),
    \label{eq:fhn-exact}
\end{equation}
where \(a\) is a parameter that sets the wave speed, \(\theta\) is the resulting propagation speed, and \(x_c\) is the initial location of the \(u=0.5\) level. For \(\theta>0\) (i.e., \(a<0.5\)), the front is an increasing sigmoid that moves to the left with constant speed~\(\theta\).
 
The derivative \(u_x\) of the increasing profile integrates to one over the line, so the globs are given equal weights \(1/N\) and are placed at the midpoint quantiles of the profile, the points at which it takes the values \(q_i\) below.

\begin{equation}
    q_i=\frac{i-\tfrac12}{N},\qquad
    X_i^0=x_c-2\log\left(q_i^{-1}-1\right),\qquad
    w_i^0=\frac{1}{N},
    \quad i=1,\ldots,N.
    \label{eq:fhn-quantile-initialization}
\end{equation}
Here $u^{\mathrm{ex}}(X_i^0,0)=q_i$. Substituting the front~\eqref{eq:fhn-exact} into~\eqref{eq:fhn-pde} determines the reaction term, using the derivatives collected in Appendix~\ref{app:benchmarks},

\begin{equation}
    f(u)
    =u(1-u)
    \left[
      \frac{\theta}{2}
      -\frac{D(1-2u)}{4}
    \right].
    \label{eq:fhn-reaction}
\end{equation}
Differentiating~\eqref{eq:fhn-reaction} gives the weight update function~\(R(u)=f'(u)\) used in~\eqref{eq:reaction-weight-update}.
\begin{equation}
    R(u)=f'(u)
    =-\frac{3D}{2}u^2
     +\left(\frac{3D}{2}-\theta\right)u
     +\left(\frac{\theta}{2}-\frac{D}{4}\right).
    \label{eq:fhn-statistic}
\end{equation}
For the parameter values used in this study (\(a=0.25\), \(D=0.5\)), the cubic~\eqref{eq:fhn-reaction} factors as
\[
    f(u)=\tfrac{D}{2}\,u(1-u)(u-u_\star),
    \qquad
    u_\star=\tfrac12-\frac{\theta}{D}\approx-0.21,
\]
with \(u_\star\) lying outside \([0,1]\). The reaction is therefore positive throughout \(0<u<1\) (a monostable nonlinearity, with \(u=0\) unstable and \(u=1\) stable). The classical bistable Nagumo cubic has an unstable third root inside \((0,1)\)~\citep{mascagni1995}.
 
The function~\(R(u)\) has the zero-integral property
\[
    \int_0^1R(u)\,du=f(1)-f(0)=0.
\]
For the whole-line front connecting the equilibrium states $0$ and $1$, this gives zero net reaction contribution to the total gradient weight. On a finite interval, that contribution is $f(u(L,t))-f(u(0,t))$. In the discrete particle system, deviations from this balance arise from approximating the integral with a finite number of particles and from the time-stepping discretization. When homogeneous Neumann reflection is active, boundary reflections can also modify the discrete total because the weight sign changes upon reflection. The resulting weights are carried directly into the next time step.

\subsection{Burgers' equation through the Cole--Hopf transformation}
\label{sec:burgers-benchmark}
 
The viscous Burgers equation
\begin{equation}
    u_t+u u_x=\nu u_{xx}
    \label{eq:burgers-pde}
\end{equation}
admits a time-independent (stationary) solution of the form
\begin{equation}
    u^{\mathrm{ex}}(x)
    =-A\tanh
    \left(
      \frac{A(x-x_c)}{2\nu}
    \right),
    \label{eq:burgers-shock}
\end{equation}
where \(A\) is the magnitude of the asymptotic states, \(x_c\) is the shock center, and \(2\nu/A\) is the characteristic width of the transition. Throughout this paper, the term \emph{shock} refers to this viscous profile connecting the left state \(u=A\) to the right state \(u=-A\). For \(\nu>0\) the profile is smooth. In the inviscid limit \(\nu\to0\), the profile approaches a discontinuity. As emphasized by Aref and Daripa~\citep{arefdaripa1984}, Burgers' equation is a valuable numerical test case since it is directly solvable through the Cole--Hopf transformation, so nonlinear discretization errors can be measured against an exact solution. Because the far-field states are equal and opposite, the profile is stationary. Appendix~\ref{app:benchmarks} verifies that~\eqref{eq:burgers-shock} satisfies~\eqref{eq:burgers-pde} with zero propagation speed.

One possible direct GRW formulation of Burgers' equation represents its nonlinear term by a weight update for the gradient~\(v=u_x\). Differentiating~\eqref{eq:burgers-pde} with respect to~\(x\) gives
\[
    v_t=\nu v_{xx}-(uv)_x=\nu v_{xx}-u_x v-u v_x,
\]
and writing the non-diffusive part as \(Rv\) with \(v_x=u_{xx}\) yields the weight update function
\[
    R=-\left(\frac{u\,u_{xx}}{u_x}+u_x\right).
\]
The ratio~\(u_{xx}/u_x\) appearing in this expression is poorly conditioned. In particular, away from the shock center, both~\(u_x\) and~\(u_{xx}\) approach zero, so a small absolute error in either numerically computed derivative can produce a large error in the ratio and hence in~\(R\). This sensitivity motivates the transformed formulation used here. A direct random-walk treatment of Burgers' equation is nevertheless possible when the nonlinear term is handled by particle advection rather than by weight updates. Roberts~\citep{roberts1989} proved convergence of such a fractional-step method. In this study, we apply GRW to the heat equation obtained through the Cole--Hopf transformation.

The Cole--Hopf transformation, due independently to Hopf~\citep{hopf1950} and Cole~\citep{cole1951}, introduces a new variable~\(\phi\) defined by
\begin{equation}
    u=-2\nu\frac{\phi_x}{\phi}.
    \label{eq:cole-hopf-recovery}
\end{equation}
Under this transformation, \(\phi\) satisfies the heat equation
\begin{equation}
    \phi_t=\nu\phi_{xx}.
    \label{eq:cole-hopf-heat}
\end{equation}
We apply the standard GRW diffusion method to evolve~\(\phi_x\)-globs for~\eqref{eq:cole-hopf-heat}, and then recover the physical Burgers field~\(u\) from the reconstructed~\(\phi\) via~\eqref{eq:cole-hopf-recovery}.
 
Given an initial Burgers field~\(u_0\), the transformed initial condition is constructed by numerical integration as
\begin{equation}
    \Psi_0(x)=\int_0^x u_0(s)\,ds,
    \qquad
    \phi_0(x)
    =C\exp\left[-\frac{\Psi_0(x)}{2\nu}\right],
    \label{eq:cole-hopf-initialization}
\end{equation}
where \(\Psi_0\) is evaluated by the trapezoidal rule on \(P\) initialization points~\(\{x_j\}\), and \(C>0\) is an arbitrary normalization constant that cancels from~\eqref{eq:cole-hopf-recovery}. In the implementation, we set \(C\) so that \(\max_x\phi_0(x)=1\).

For the centered stationary shock, $k=A/(2\nu)$ and $x_c=L/2$ give
\[
\phi_0(x)=\frac{\cosh(k(x-x_c))}{\cosh(kx_c)},\qquad
\phi(x,t)=e^{\nu k^2t}\phi_0(x).
\]
The spatially uniform prefactor cancels in the recovery ratio, so the Burgers profile remains stationary while the transformed endpoint values grow. This distinction is tested by the boundary control of Section~\ref{sec:colehopf-boundary}. The derivation is given in Appendix~\ref{app:benchmarks}. At $A=1$, $\nu=0.5$, $T=0.5$, and $L=4$, the minimum of $\phi$ at the shock center is $1/\cosh2=0.266$ initially and $e^{0.25}/\cosh2=0.341$ at the final time.

We then place the \(\phi_x\)-globs at the midpoints of adjacent initialization points with forward-difference weights
\begin{equation}
    w_j^\phi
    =\phi_0(x_{j+1})-\phi_0(x_j),
    \label{eq:phi-gradient-weight}
\end{equation}
which sum to~\(\phi_0(L)-\phi_0(0)\) because the intermediate terms cancel. This placement yields one glob per interval, so the \(P\) initialization points produce \(P-1\) \(\phi_x\)-globs. These globs are evolved by~\eqref{eq:brownian-update} with \(\kappa=\nu\).
 
At the final time, we use $M$ output locations $x_j=jh$, $h=L/(M-1)$, and assign each transported weight to index $\min(\lfloor X_i/h\rfloor,M-1)$. These indexed bins follow the output-point spacing, while the heat bin-and-sum studies use $M$ equal intervals of width $L/M$. The binned weights are smoothed with a Gaussian kernel whose standard deviation is fixed at 12 reconstruction-point spacings, so its physical width is 
\[
\sigma_x=12\frac{L}{M-1}.
\] 
Near the domain boundaries, where the kernel is truncated by the domain edge, the smoothed values are divided by the kernel mass remaining inside the domain, so that the truncated kernel integrates to one again to compensate for the missing kernel mass. This normalization alone does not eliminate all boundary smoothing error. This width suppresses stochastic particle noise while remaining narrow relative to the shock transition width~\(2\nu/A\) at the representative configuration, where $\sigma_x\approx0.120$ at $M=400$. Section~\ref{sec:colehopf-decoupled} varies this bandwidth independently. The smoothed weights are then corrected so that their sum equals the prescribed endpoint difference~\(\phi_0(L)-\phi_0(0)\). When this difference is zero, as for the stationary shock, the mean of the smoothed weights is subtracted, because no rescaling maps a nonzero sum to zero. Otherwise, the weights are rescaled proportionally. If the corrected bin weights are~\(\widehat w_k^\phi\), then the cumulative reconstruction of the transformed field is
\begin{equation}
    \phi_N(x_j)
    =\phi_0(0)+\sum_{k=0}^{j}\widehat w_k^\phi.
    \label{eq:phi-reconstruction}
\end{equation}
The numerical derivative~\((\phi_N)_x\) is computed from the smoothed reconstruction, and the Burgers field is recovered as
\[
u_N=-2\nu\frac{(\phi_N)_x}{\phi_N},
\]
the discrete form of~\eqref{eq:cole-hopf-recovery}. A small positive floor is applied to~\(\phi_N\) in the denominator to prevent division by zero where particle noise drives the reconstruction near or below zero. The floor is set to half the minimum value of the initial transformed profile~\(\phi_0\), with an absolute lower bound of \(10^{-10}\). For the canonical shock, the closed form above gives a floor of \(1/(2\cosh 2)=0.133\). The exact field starts a factor of two above this floor by construction and its minimum grows by the prefactor \(e^{\nu k^{2}t}\), so activation indicates an error in the reconstructed transformed field. Wherever this floor is active, the recovered field is set to zero. Because the recovery step differentiates the reconstructed~\(\phi\) and divides by~\(\phi_N\), small errors in the reconstructed transformed field are amplified, especially in regions where~\(\phi_N\) is small.
 
Algorithm~\ref{alg:cole-hopf} summarizes the complete Cole--Hopf GRW procedure. The GRW solver operates on the transformed field~\(\phi\). The nonlinear Burgers field~\(u\) appears only at the initialization and recovery stages.

\FloatBarrier
\begin{algorithm}[!htbp]
\caption{Cole--Hopf GRW procedure for Burgers' equation~\eqref{eq:burgers-pde}}
\label{alg:cole-hopf}
\KwIn{viscosity \(\nu\), initial field \(u_0\), time step \(\Delta t\), number of steps \(K\), \(P\) initialization points, \(M\) reconstruction bins.}
\vspace{0.15em}
\KwOut{recovered Burgers field on the \(M\) reconstruction bins.}
\vspace{0.3em}
Compute \(\Psi_0\) and \(\phi_0\) on the \(P\) initialization points via Eq~\eqref{eq:cole-hopf-initialization}. \\
Place the \(P-1\) \(\phi_x\)-globs at the midpoints of adjacent initialization points with weights via Eq~\eqref{eq:phi-gradient-weight}. \\
\vspace{0.3em}
\For{\(n = 0, \ldots, K-1\)}{
    \textit{Diffuse:} displace each glob via Eq~\eqref{eq:brownian-update} with \(\kappa=\nu\). \\
    \textit{Reflect:} reflect any glob that exits \([0,L]\), leaving the weight unchanged. \\
}
\vspace{0.3em}
Assign each weight to bin \(k=\min(\lfloor X_i/h\rfloor,\,M-1)\). \\
Smooth with a Gaussian kernel of standard deviation 12 reconstruction-point spacings. Near the boundaries, divide the smoothed values by the kernel mass remaining inside the domain to correct boundary bias. \\
Correct the smoothed weights so that their sum equals \(\phi_0(L)-\phi_0(0)\). When this endpoint difference is zero to within numerical tolerance, subtract the mean. Otherwise, rescale the weights proportionally. \\
Reconstruct \(\phi_N\) via Eq~\eqref{eq:phi-reconstruction}. \\
Compute the numerical derivative \((\phi_N)_x\) from the smoothed reconstruction. \\
Recover the field via Eq~\eqref{eq:cole-hopf-recovery}, flooring the denominator at \(\max(\tfrac12\min_x\phi_0,\,10^{-10})\) and setting the recovered value to zero wherever the floor is active. \\
\end{algorithm}
\FloatBarrier

\section{Verification and Ensemble Methodology}
\label{sec:verification}

This section fixes the measurement conventions used throughout. All studies share the discrete norms defined first. The representative calculations use the single fixed seed 42, which makes each of them exactly repeatable. Every multi-seed study is an ensemble over a fixed seed list, while the deterministic controls use no random input. The following subsections define the ensemble vocabulary, the bias--spread decomposition, the fitted rates, and their uncertainty.

\subsection{Discrete profile errors}
\label{sec:norms}
 
All particle outputs are first reconstructed at the finite uniform reconstruction points introduced in Section~\ref{sec:reconstruction}, with \(M\)~points \(\{x_j\}_{j=0}^{M-1}\) and spacing~\(h\). For the pointwise error
\[
    e_j=u_N(x_j)-u^{\mathrm{ref}}(x_j),
\]

we report the discrete norm
\vspace{0.6em}
\begin{equation}
    \lVert e\rVert_{L_h^2}
    =\left(
      h\sum_{j=0}^{M-1}|e_j|^2
    \right)^{1/2}.
    \label{eq:grid-l2}
\end{equation}
\vspace{0.6em}

The spacing factor~\(h\) makes this conventional grid norm approximate the continuous $L^2$ norm~\citep{leveque2007}. When we need an average pointwise error scale, we use the root mean square error
\vspace{0.6em}
\begin{equation}
    \operatorname{RMSE}(e)
    =\left(
      \frac{1}{M}
      \sum_{j=0}^{M-1}|e_j|^2
    \right)^{1/2}.
    \label{eq:rmse}
\end{equation}
\vspace{0.6em}

For Burgers, the domain, boundary, and decoupled studies report RMSE, the coupled plateau and amplification response report $L_h^2$, and the inversion control gives both. Tables and figures identify the norm used. The two are related exactly by
\[
\lVert e\rVert_{L_h^2}=\sqrt{hM}\,\operatorname{RMSE}(e).
\]
For uniform endpoint-inclusive reconstruction points, including the Burgers output, $h=L/(M-1)$. For the heat bin-and-sum studies, $h=L/M$, so in both cases $\sqrt{hM}\approx\sqrt L$, which is $2$ for the $L=4$ configurations used below.

\subsection{Realizations, seeds, and ensembles}
\label{sec:ensemble-vocab}

A realization is a complete run driven by one pseudorandom seed. An ensemble comprises $S$ runs with identical numerical parameters and different seeds, treated as independent realizations. Each study reuses the same list of seeds at every particle count $N$, so realizations at different counts reuse pseudorandom number sequences, with the assignment of increments to globs depending on $N$. Section~\ref{sec:ensemble-bootstrap} states how the bootstrap treats this sharing.

\subsection{The bias--spread decomposition and its identity}
\label{sec:ensemble-errors}

For a fixed particle count $N$, let $u_{N,s}$ denote the numerical field produced by realization $s$, evaluated at the study's reconstruction points. For an ensemble of $S$ realizations, define the mean field
\[
\bar u_N(x)
=
\frac{1}{S}
\sum_{s=1}^{S}u_{N,s}(x),
\]
i.e., the pointwise average of the $S$ reconstructed profiles. Relative to a reference solution $u^{\mathrm{ref}}$, we then define three error quantities.
\begin{subequations}
\label{eq:error-decomp}
\begin{align}
E_{\mathrm{bias}}(N)
&=
\|\bar u_N-u^{\mathrm{ref}}\|_{L_h^2},
\label{eq:E-bias}\\
E_{\mathrm{spread}}(N)
&=
\left[
\frac{1}{S}
\sum_{s=1}^{S}
\|u_{N,s}-\bar u_N\|_{L_h^2}^{2}
\right]^{1/2},
\label{eq:E-spread}\\
E_{\mathrm{total}}(N)
&=
\left[
\frac{1}{S}
\sum_{s=1}^{S}
\|u_{N,s}-u^{\mathrm{ref}}\|_{L_h^2}^{2}
\right]^{1/2}.
\label{eq:E-total}
\end{align}
\end{subequations}
The ensemble-mean error $E_{\mathrm{bias}}$ estimates systematic discrepancy but retains finite-ensemble sampling error. The spread $E_{\mathrm{spread}}$ measures variation about the mean, and $E_{\mathrm{total}}$ is the root-mean-square realization error. For independent, identically distributed realizations with $S\ge2$, let $b=\mathbb E[u_{N,s}]-u^{\mathrm{ref}}$ and $V=\mathbb E\|u_{N,s}-\mathbb E[u_{N,s}]\|_{L_h^2}^2$. Then
\[
\mathbb E[E_{\mathrm{spread}}^2]=\frac{S-1}{S}V,\qquad
\mathbb E[E_{\mathrm{bias}}^2]=\|b\|_{L_h^2}^2+\frac VS
=\|b\|_{L_h^2}^2+\frac{\mathbb E[E_{\mathrm{spread}}^2]}{S-1}.
\]
The second identity shows that $E_{\mathrm{bias}}$ contains a sampling contribution $V/S$ in addition to the systematic error $\|b\|_{L_h^2}$. We report the descriptive statistics~\eqref{eq:error-decomp} and use these identities to interpret them. Appendix~\ref{app:heat-uncertainty} derives the finite-ensemble factors and the uncertainty calculations.

Monte Carlo sampling provides the $N^{-1/2}$ comparison rate for the stochastic spread~\citep{caflisch1998}. The bias has no universal rate because it depends on systematic effects introduced by the numerical procedure. Measuring only the total error would combine these two distinct mechanisms.

The three quantities are not independent. They satisfy the discrete bias--spread identity
\begin{equation}
E_{\mathrm{total}}^2
=
E_{\mathrm{bias}}^2
+
E_{\mathrm{spread}}^2,
\label{eq:bias-variance}
\end{equation}
up to floating-point arithmetic. We also use this identity to check the heat bias--spread--total statistics reported in Sections~\ref{sec:heat-ensemble} and~\ref{sec:heat-grid-paired}.

Indeed, $u_{N,s}-u^{\mathrm{ref}}=(u_{N,s}-\bar u_N)+(\bar u_N-u^{\mathrm{ref}})$, and averaging its squared norm eliminates the cross term because $\sum_s(u_{N,s}-\bar u_N)=0$. Equation~\eqref{eq:bias-variance} is therefore an identity for the realized statistics, distinct from the expectation identities above.

\subsection{Empirical convergence rates}
\label{sec:ensemble-rates}

We summarize particle-count dependence by the power law $E(N)\approx cN^r$. Taking logarithms gives
\begin{equation}
\log E(N)
=
\log c+r\log N.
\label{eq:loglog-rate-fit}
\end{equation}
We estimate $r$ by least squares on $(\log N_k,\log E(N_k))$ over the tested counts. We call $r$ the \emph{fitted exponent} or \emph{observed rate} throughout this paper.

The reference value is the Monte Carlo rate $r=-1/2$, the rate at which the statistical error of an average of $N$ independent samples decreases~\citep{caflisch1998}. For the heat problem, independence makes the squared sampling contribution proportional to $1/N$ at fixed comparison points, so the root-mean-square spread scales as $N^{-1/2}$. Section~\ref{sec:analysis-cdf} derives this relation. In the reaction solver the globs interact through the reconstructed field, so there $-1/2$ serves as a comparison value. The same theoretical rate appears in error estimates for gradient-based Monte Carlo methods in the hyperbolic relaxation setting~\citep{bertaglia2024}. These fitted exponents describe behavior over the stated finite refinement ranges.

\subsection{Uncertainty of the fitted slopes}
\label{sec:ensemble-bootstrap}

Uncertainty in the fitted trends is quantified using a realization-level percentile bootstrap~\citep{efron1981}. Each study reuses one seed list across its parameter values. We resample the seed labels with replacement and use the same resampled labels at every value, preserving the grouping of realizations across the study. For each of $5000$ bootstrap replicates, we recompute the reported statistic and refit the trend, using the log--log slope for the particle-count and initialization-point rate studies and an ordinary linear fit for the Burgers domain trend against $L$. The heat and scalar reaction--diffusion rate studies use thirty realizations per particle count, the Burgers decoupled study twenty, and the Burgers domain study thirty. This construction propagates seed-to-seed variability without assuming a normal distribution for the fitted trend. Appendix~\ref{app:bootstrap} compares these joint intervals with intervals obtained by resampling independently within each parameter value.

\subsection{Reaction-front location diagnostics}
\label{sec:front-diagnostics}
 
For the traveling front~\eqref{eq:fhn-exact}, we define the front location $x_{\mathrm{center}}$ as the point where the field reaches the midpoint value:
\begin{equation}
    u(x_{\mathrm{center}}(t),t)=0.5,
    \qquad
    x_{\mathrm{center}}^{\mathrm{exact}}(t)=x_c-\theta t.
    \label{eq:fhn-front-location}
\end{equation}
The exact front location moves to the left at constant speed~\(\theta\). For the multi-seed convergence study of Section~\ref{sec:fhn-convergence}, the numerical front location is the linearly interpolated position at which the reconstructed profile crosses \(u=1/2\). If several crossings are detected, the middle crossing is selected. If no crossing is present, the reconstruction point nearest to \(u=1/2\) is used. In the plot of front location versus time in Section~\ref{sec:fhn-representative}, the crossing is instead extracted directly from the sorted glob data, as the position of the first glob at which the cumulative weight sum reaches~\(0.5\). There we set the reference at the initial crossing of the discrete reconstruction at~\(t=0\), so that the plotted displacement measures the front motion relative to the particle profile at the initial time.
 
This diagnostic is inherently local. A small change in the placement of a few globs near the~\(u=0.5\) level, or a change in the interpolation, can shift the extracted location even when the overall shape of the profile is unchanged. Front-location error is therefore expected to be noisier than the global~\(L_h^2\) profile error, and we treat it as a secondary measure.

\subsection{Cole--Hopf Burgers error decomposition}
\label{sec:burgers-decomposition}
 
For the Burgers problem, the total error in the recovered field~\(u\) accumulates across multiple stages of the Cole--Hopf procedure. To separate these, we define three RMSE quantities. Let \(u^{\mathrm{ex}}\) denote the exact stationary shock~\eqref{eq:burgers-shock}, let \(u^{\mathrm{FD}}\) denote the field recovered from a deterministic (i.e., non-stochastic) finite-domain solution of~\eqref{eq:cole-hopf-heat} using the same endpoint values as the GRW calculation, and let \(u^{\mathrm{GRW}}\) denote the field recovered from the GRW particle calculation.

The deterministic reference solves~\eqref{eq:cole-hopf-heat} by forward Euler in time and second-order central differences in space~\citep{leveque2007}, on the same reconstruction points with the endpoint values of~\(\phi\) held fixed at every step. Its internal time step is set by the stability constraint \(\nu\,\Delta t_{\mathrm{FD}}/h^2\leq0.4\) rather than by the GRW time step. The two computations share these points and the fixed endpoint data but differ in their temporal discretization and discrete endpoint enforcement.

The GRW reconstruction uses~\(\phi_0(0)\) as the integration constant and enforces the prescribed endpoint difference through the corrected weight sum. This reconstruction difference is included in~\(E_{\mathrm{GRW}}\). The derivative and the recovery~\eqref{eq:cole-hopf-recovery} are applied to the deterministic~\(\phi\) with the same central differences used for the GRW reconstruction. We then define
\begin{subequations}
\label{eq:burgers-error-decomposition}
\begin{align}
    E_{\mathrm{det}}
    &=
    \operatorname{RMSE}
    \left(u^{\mathrm{FD}}-u^{\mathrm{ex}}\right),
    \label{eq:domain-mismatch}\\
    E_{\mathrm{GRW}}
    &=
    \operatorname{RMSE}
    \left(u^{\mathrm{GRW}}-u^{\mathrm{FD}}\right),
    \label{eq:grw-reconstruction-error}\\
    E_{\mathrm{total}}
    &=
    \operatorname{RMSE}
    \left(u^{\mathrm{GRW}}-u^{\mathrm{ex}}\right).
    \label{eq:recovered-field-error}
\end{align}
\end{subequations}
Here \(E_{\mathrm{det}}\) is the error that remains even when the transformed heat equation is solved without any particle noise. Specifically, it includes the mismatch between the fixed transformed endpoint data and the exact shock, together with deterministic discretization and recovery errors. The quantity~\(E_{\mathrm{GRW}}\), which we call the \emph{particle error}, measures the additional error introduced by replacing the deterministic computation with the GRW particle calculation. This includes stochastic particle noise, binning, Gaussian smoothing, and the correction of the smoothed weights to enforce the prescribed endpoint difference. The quantity~\(E_{\mathrm{total}}\) is the total error in the recovered Burgers field.
 
These quantities compare the deterministic and particle effects at a given parameter set. The corresponding error fields add to the total error field, but their RMSE values do not generally add to~\(E_{\mathrm{total}}\) because the fields are not orthogonal. We therefore use~\(E_{\mathrm{det}}\) and~\(E_{\mathrm{GRW}}\) to identify the larger contribution.

\section{Heat Equation Results}
\label{sec:heat}

The heat study first examines reconstruction and sampling errors, then uses paired comparisons to identify the fixed-bin floor. The final study uses the resulting identity to predict errors and select a particle count for a second initial profile.

\subsection{Representative reconstruction}
\label{sec:heat-representative}

\begin{figure}[pos=htbp]
    \centering
\includegraphics[width=0.85\textwidth]{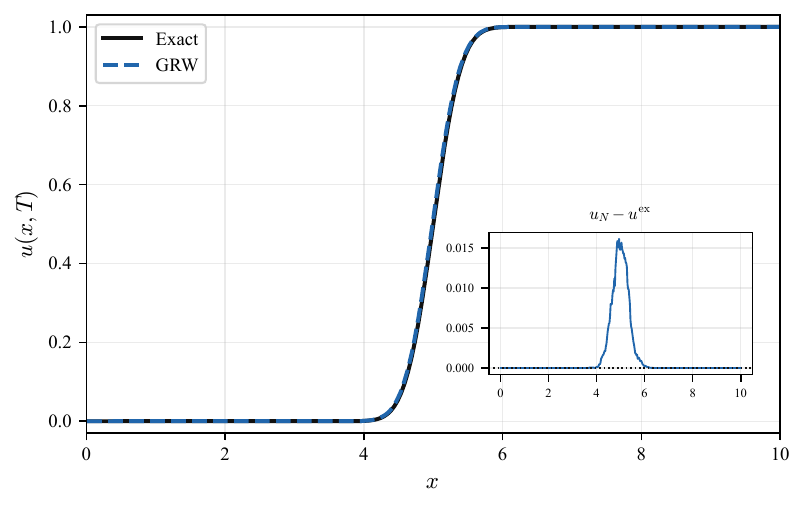}
    \caption{Heat equation. GRW reconstruction versus exact error-function profile for \(N=50{,}000\), \(\alpha=0.1\), \(T=0.5\). The GRW reconstruction (dashed line) lies nearly on top of the exact profile, and the two curves are visually indistinguishable at this particle count. The inset shows the residual \(u_N-u^{\mathrm{ex}}\), which is concentrated near the transition region.}
    \label{fig:heat-comparison}
\end{figure}

We first provide results for the heat-equation step problem of Section~\ref{sec:heat-benchmark}, with \(\alpha=0.1\), domain~\([0,10]\), step location~\(x_0=5\), step values \(u_L=0\) and \(u_R=1\), final time~\(T=0.5\), time step~\(\Delta t=0.001\), and \(N=50{,}000\)~globs. Both boundaries use the weight-preserving (Dirichlet-type) reflection, under which the reconstructed field keeps the values \(u=0\) at the left boundary and \(u=1\) at the right boundary. At the final time, the diffusion length is \(\sqrt{2\alpha T}\approx0.32\), which is small compared with the distance of~5 from the step to either boundary, so the finite computational domain has a negligible effect on the comparison with the infinite-domain solution~\eqref{eq:heat-exact}. The reconstruction for Figure~\ref{fig:heat-comparison} uses 400 reconstruction bins.
 
\FloatBarrier
 
The residual in Figure~\ref{fig:heat-comparison} is concentrated near the transition, where the gradient is largest. A comparison of two reconstructions of this same particle set supplies the first diagnostic of the error's origin.

Reconstructing this same $N=50{,}000$ run on 300 reconstruction bins instead of 400 raises the error from $L_h^2=0.0119$ to $0.0158$. The ratio $0.0119/0.0158\approx0.75$ matches the ratio of the bin widths, $300/400$, and provides a first single-seed indication of a reconstruction-bin effect at large particle counts. Section~\ref{sec:heat-grid-paired} tests this effect under a paired ensemble design and identifies its source in the evaluation convention of the bin-and-sum reconstruction.

\subsection{Multi-seed ensemble under coupled refinement}
\label{sec:heat-ensemble}

The ensemble study uses the increasing unit step $u_L=0$, $u_R=1$ on $[0,4]$, with
\begin{equation}
    \alpha=0.5,
    \qquad
    x_0=2,
    \qquad
    T=0.5,
    \qquad
    \Delta t=0.005.
    \label{eq:heat-study-parameters}
\end{equation}
These are the diffusion coefficient, domain, final time, and time step of the transformed Burgers calculation in Section~\ref{sec:burgers-representative}, so the ensemble provides a heat-only reference for sampling and reconstruction at the same numerical scales.

The reflected Gaussian steps sample the exact final-time particle distribution, as explained in Section~\ref{sec:boundary-reflection}. The reference comparison therefore measures sampling, reconstruction, and domain effects. The reference~\eqref{eq:heat-exact} is the solution on the whole line, while the computation reflects particles at the walls of $[0,4]$. We call the resulting difference the \emph{domain truncation error}. At this final time the diffusion length $\sqrt{2\alpha T}\approx0.71$ places the walls about $2.8$ standard deviations from the step, so this error is small, and the deterministic control of Section~\ref{sec:heat-grid-paired} computes it exactly as $1.1\times10^{-3}$.

Thirty stochastic realizations were run for each particle count in
\begin{equation}
    N\in
    \left\{
        500,\,
        1000,\,
        2000,\,
        5000,\,
        10{,}000,\,
        20{,}000,\,
        50{,}000
    \right\},
    \label{eq:heat-N-sequence}
\end{equation}
with the same sequence of thirty seed integers reused at every particle count, as defined in Section~\ref{sec:ensemble-vocab}. Each realization is reconstructed and compared with the exact profile~\eqref{eq:heat-exact} at a uniform set of $N$ reconstruction points on $[0,4]$. The reconstruction points are common to all realizations at each count, with $M=N$. We call this equality of reconstruction-point and particle counts \emph{coupled refinement}. Section~\ref{sec:heat-grid-paired} contrasts it with reconstructions whose bin count is held fixed. Table~\ref{tab:heat-convergence} reports the resulting bias--spread--total decomposition of Section~\ref{sec:ensemble-errors}, with each entry computed by evaluating~\eqref{eq:error-decomp} over the thirty stored realization profiles.

\begin{table}[pos=htbp]
    \centering
    \small
    \caption{Heat GRW bias--spread--total decomposition for the increasing unit-step problem ($\alpha=0.5$, $T=0.5$, $\Delta t=0.005$, $S=30$). The spread dominates the bias at every particle count.}
    \label{tab:heat-convergence}
    \begin{tabular}{rccc}
        \toprule
        $N$
        & $E_{\mathrm{bias}}$
        & $E_{\mathrm{spread}}$
        & $E_{\mathrm{total}}$\\
        \midrule
        500   & 0.00530 & 0.02850 & 0.02899\\
        1000  & 0.00261 & 0.02040 & 0.02057\\
        2000  & 0.00211 & 0.01393 & 0.01409\\
        5000  & 0.00194 & 0.00880 & 0.00901\\
        10000 & 0.00220 & 0.00703 & 0.00736\\
        20000 & 0.00156 & 0.00486 & 0.00511\\
        50000 & 0.00113 & 0.00330 & 0.00349\\
        \bottomrule
    \end{tabular}
\end{table}

\FloatBarrier

\begin{figure}[pos=htbp]
    \centering
\includegraphics[width=0.90\textwidth]{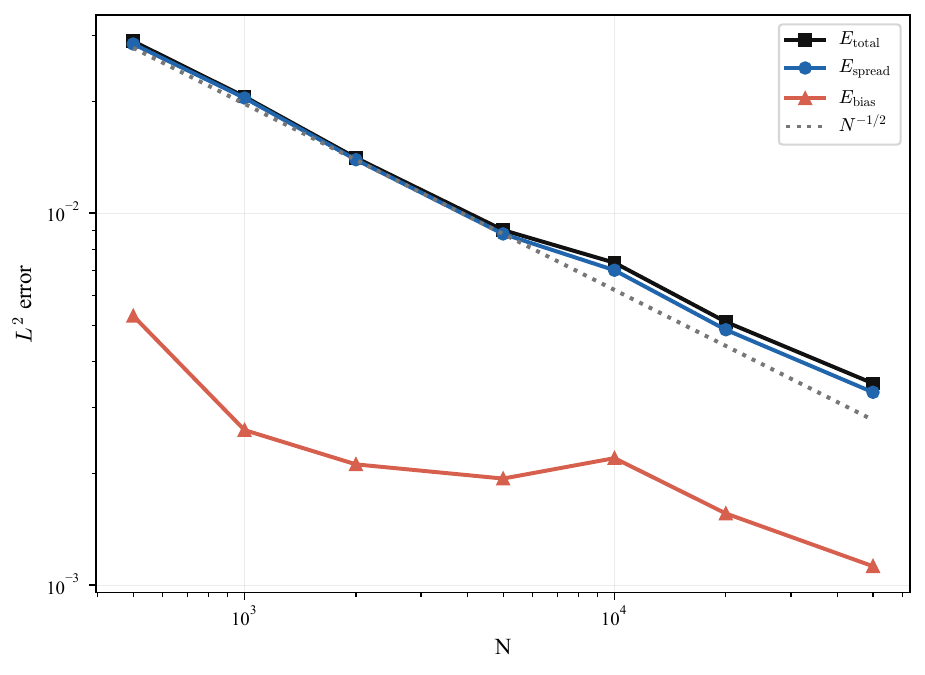}
    \caption{Bias, stochastic spread, and total error for the heat GRW refinement study on log--log axes. The spread is the dominant contribution and has fitted exponent $-0.468$. The total error has fitted exponent $-0.457$.}
    \label{fig:heat-bst}
\end{figure}

Reading down the columns of Table~\ref{tab:heat-convergence}, the spread falls steadily by a factor of $8.6$ across the hundredfold increase in $N$. The bias is between $2.9$ and $7.8$ times smaller than the spread across the tested counts, and the total error tracks the spread closely. This is exactly the pattern expected when sampling noise dominates. A least-squares fit of $\log E$ against $\log N$, following Section~\ref{sec:ensemble-rates}, gives the fitted exponents
\begin{subequations}
\label{eq:heat-convergence-rates}
\begin{align}
    E_{\mathrm{total}}&\sim N^{-0.457},
    \label{eq:heat-total-rate}\\
    E_{\mathrm{spread}}&\sim N^{-0.468}.
    \label{eq:heat-spread-rate}
\end{align}
\end{subequations}
Following Section~\ref{sec:ensemble-bootstrap}, the $95\%$ confidence intervals on the fitted exponents~\eqref{eq:heat-convergence-rates} are $[-0.521,-0.416]$ for the spread and $[-0.505,-0.410]$ for the total error. Both contain the Monte Carlo rate $-1/2$. The bias remains below the spread at every count and has a less regular trend. Figure~\ref{fig:heat-bst} displays the decomposition graphically.

\FloatBarrier

\subsection{Paired reconstruction study of the fixed-bin floor}
\label{sec:heat-grid-paired}

The coupled ensemble above evaluates the cumulative reconstruction directly at $M=N$ uniform points, refining these points together with the particle count, while the representative diagnostic of Section~\ref{sec:heat-representative} shows a floor when the field is instead binned on a fixed set of bins. To attribute that floor, we repeat the ensemble of Section~\ref{sec:heat-ensemble} as a paired experiment. For each particle count and each of the same thirty seeds, one simulation is run, and the identical final particle set is then reconstructed under five treatments (a treatment being one reconstruction and comparison convention applied to the same particles).

One is the coupled treatment, direct evaluation at $M=N$ points as in Section~\ref{sec:heat-ensemble}. Two apply the bin-and-sum treatment of the representative study at $300$ and $400$ fixed bins, with the exact profile evaluated at the bin centers (the convention of Section~\ref{sec:heat-representative}). The remaining two repeat those fixed-bin reconstructions with the exact profile evaluated at the bin right edges instead.

The last convention is motivated by what the reconstruction actually is. The cumulative sum over bins $0,\dots,k$ counts every glob up to the right edge of bin $k$, so it is the reconstruction at that edge, and comparing it at the bin center misplaces the evaluation by half a bin. The five treatments use identical trajectories for each seed, so each comparison changes only the reconstruction or evaluation convention.

A deterministic bin-and-sum control accompanies the ensemble. Under the weight-preserving (Dirichlet-type) reflection on $[0,L]$, the distribution function $F_R$ of a particle position at time $T$, introduced in Section~\ref{sec:boundary-reflection}, is given by the method of images as a sum of Gaussian images reflected across both walls~\citep{carslawjaeger1959}. Appendix~\ref{app:heat-uncertainty} states the density. We evaluate $F_R$ through the truncated antiderivative
\[
G(x)=\sum_{n=-4}^{4}\ \sum_{\mu\in\{2nL+x_0,\ 2nL-x_0\}}\frac12\left[1+\operatorname{erf}\!\left(\frac{x-\mu}{\sqrt{4\alpha T}}\right)\right],
\]
where $\mu$ runs over the images of $x_0$ under repeated reflection across the two walls, with the sum truncated at $\lvert n\rvert\le4$, beyond which the images lie many diffusion lengths outside $[0,L]$ and contribute below machine precision. The truncated sum $G$ is an antiderivative rather than a normalized distribution function. On $[0,L]$, $F_R(x)=G(x)-G(0)$, with $F_R(L)=1$ to machine precision at this truncation. The bin masses below are differences of $G$, so the additive constant cancels. Binning $F_R$ gives the exact expected mass $G(x_{k+1})-G(x_k)$ in bin $k$, and the cumulative sum of those masses is the deterministic reconstruction, formed exactly as in the particle treatments but with no random input. Comparing it with the exact profile predicts
\[
\text{bin-center}\ \ 0.00435\ (300),\ \ 0.00334\ (400),\qquad \text{bin-edge}\ \ 0.00109\ (\text{both}).
\]
The bin-edge value equals, to three digits, the norm of the difference between the reflected distribution $F_R$ on $[0,L]$ and the infinite-line profile~\eqref{eq:heat-exact}, computed on a fine grid. This is the domain truncation error defined in Section~\ref{sec:heat-ensemble}. The bin-center floor is first order in the bin width because the cumulative sum is the reconstruction at the right edge while the reference is evaluated half a bin to its left, which by the Taylor expansion~\eqref{eq:halfbin} of Section~\ref{sec:analysis-coordinates} costs about $(h/2)\lvert u_x\rvert$ at each point. The predicted $300$-to-$400$ ratio $1.30$ is close to the bin-width ratio $4/3$, the departure reflecting the domain truncation error and its cross term with the alignment shift. For the bin-edge comparison the control predicts a value independent of the bin count.

\begin{table}[pos=htbp]
\centering
\small
\caption{Paired heat reconstruction study using the same thirty realizations per particle count. Entries are $L_h^2$ errors for the coupled points and the 300-bin reconstructions. The 400-bin results show the same pattern and are reported in the text. The spread differs by about $0.3\%$ among treatments.}
\label{tab:heat-grid-paired}
\begin{tabular}{rccccc}
\toprule
& \multicolumn{1}{c}{coupled ($M=N$)} & \multicolumn{2}{c}{fixed $300$, bin center} & \multicolumn{2}{c}{fixed $300$, right edge}\\
\cmidrule(lr){2-2}\cmidrule(lr){3-4}\cmidrule(lr){5-6}
$N$ & $E_{\mathrm{total}}$ & $E_{\mathrm{total}}$ & $E_{\mathrm{bias}}$ & $E_{\mathrm{total}}$ & $E_{\mathrm{bias}}$\\
\midrule
500   & 0.02899 & 0.02960 & 0.00768 & 0.02891 & 0.00430\\
1000  & 0.02057 & 0.02101 & 0.00495 & 0.02058 & 0.00255\\
2000  & 0.01409 & 0.01467 & 0.00466 & 0.01404 & 0.00191\\
5000  & 0.00901 & 0.01037 & 0.00548 & 0.00899 & 0.00185\\
10000 & 0.00736 & 0.00930 & 0.00608 & 0.00735 & 0.00213\\
20000 & 0.00511 & 0.00716 & 0.00526 & 0.00510 & 0.00154\\
50000 & 0.00349 & 0.00564 & 0.00457 & 0.00349 & 0.00112\\
\bottomrule
\end{tabular}
\end{table}
\FloatBarrier

With the exact profile evaluated at bin centers, the measured bias at $N=50{,}000$ is $4.57\times10^{-3}$ on 300 bins and $3.56\times10^{-3}$ on 400 bins, against the deterministic predictions of $4.35\times10^{-3}$ and $3.34\times10^{-3}$. Evaluating at the corresponding right edges reduces the measured bias to $1.12\times10^{-3}$ at both bin counts, against the predicted domain truncation error of $1.09\times10^{-3}$, below the spread throughout the tested range. The fixed-bin and coupled total errors then agree to the reported precision at $N=50{,}000$ ($0.00349$ in each case). The spread changes by about $0.3\%$ among the treatments and retains fitted exponent $-0.468$.

Because the bin-center and bin-edge comparisons evaluate the same reconstructions, their spreads are identical, and the difference of their squared biases contains the evaluation shift together with its cross term with the residual bias. Writing $d=r_{\mathrm{e}}-r_{\mathrm{c}}$ for the reference evaluated at the edges minus the reference evaluated at the centers, and $b_{\mathrm{e}}=\bar u_N-r_{\mathrm{e}}$ for the residual bias after alignment, expanding the squared norm gives
\begin{equation}
\lVert b_{\mathrm{e}}+d\rVert_{L_h^2}^2-\lVert b_{\mathrm{e}}\rVert_{L_h^2}^2=\lVert d\rVert_{L_h^2}^2+2\langle b_{\mathrm{e}},d\rangle_{h},
\label{eq:paired-identity}
\end{equation}
where $\langle\cdot,\cdot\rangle_h$ is the inner product that induces the norm~\eqref{eq:grid-l2}. At $N=50{,}000$ the shift norms are $\lVert d\rVert_{L_h^2}=4.21\times10^{-3}$ at 300 bins and $3.16\times10^{-3}$ at 400 bins, and the cross terms are $1.9\times10^{-6}$ and $1.5\times10^{-6}$, about $11\%$ and $15\%$ of $\lVert d\rVert_{L_h^2}^2$. The squared shift and its cross term therefore account for the difference between the squared measured biases, $4.43\times10^{-3}$ and $3.38\times10^{-3}$ after taking square roots, and we quantify the alignment effect by $\lVert d\rVert_{L_h^2}$ and by the deterministic control.

The measured bin-center biases exceed the deterministic predictions by $5\%$ and $7\%$. This difference is compatible with sampling fluctuation in the mean of $S=30$ realizations. The exact second-moment calculation of Appendix~\ref{app:heat-uncertainty} gives the root-mean-square values $\sqrt{\mathbb{E}[E_{\mathrm{bias}}^2]}=4.38\times10^{-3}$ at 300 bins and $3.38\times10^{-3}$ at 400 bins, and the approximate standard deviations of the measured norms are about $3.6\times10^{-4}$ and $3.5\times10^{-4}$. The measured norms lie about half a standard deviation above these values, which supports sampling fluctuation as the explanation.

\begin{figure}[pos=htbp]
    \centering
    \includegraphics[width=1\textwidth]{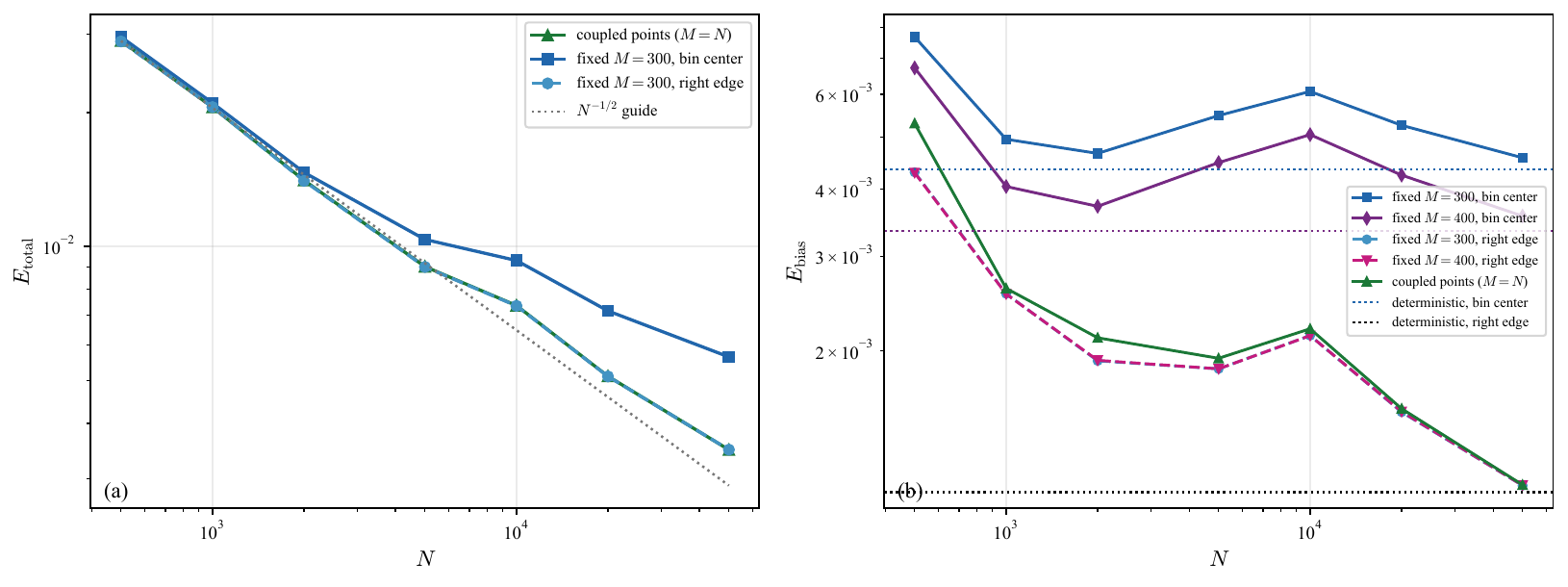}
    \caption{Paired heat reconstruction study. (a)~Total error for the coupled treatment and the fixed 300 bins when the exact profile is evaluated at bin centers or at the corresponding right edges. At the bin edges, the fixed bins track the coupled treatment. (b)~Bias component against the deterministic bin-and-sum predictions (dotted). The bin-center errors approach their predicted floors. The bin-edge errors approach the domain truncation error.}
    \label{fig:heat-grid-paired}
\end{figure}
\FloatBarrier

Figure~\ref{fig:heat-grid-paired} summarizes how the paired and deterministic controls localize the apparent fixed-bin floor to the half-bin comparison shift. Mesh-error saturation is known in gradient-based Monte Carlo refinement studies~\citep{bertaglia2024}. Here the practical issue is more specific. A cumulative reconstruction must be compared at the location it represents, or a first-order alignment error can appear as a reconstruction-bin floor.

\subsection{Predicting the floor and the particle count}
\label{sec:heat-two-step}

The controls above diagnose an error floor after a refinement study has been run. We now use an identity for the expected squared error to predict the floor and the particle count at which it begins to govern the error. This section tests those predictions on a second heat configuration, chosen so that the globs carry unequal weights and so that the two references of Section~\ref{sec:heat-ensemble} differ appreciably.

The initial condition places two steps of unequal height on $[0,L]$:
\begin{equation}
u(x,0)=u_L+A\,H(x-a)+B\,H(x-b),\qquad 0<a<b<L,
\label{eq:two-step-data}
\end{equation}
with $H$ the unit step function. We take $u_L=0$, $A=0.3$, $B=0.7$, $a=1.4$, $b=2.6$, and the parameters~\eqref{eq:heat-study-parameters} of the ensemble study, so that $L=4$ and the diffusion length $\sqrt{2\alpha T}\approx0.71$ places each wall about two diffusion lengths from the nearer step. Half the globs are released at $a$ and half at $b$, so $f_a=f_b=1/2$ and the two groups carry the unequal weights $2A/N$ and $2B/N$.

Two references are available, and they differ in what they measure. The finite-interval solution
\begin{equation}
m(x,T)=u_L+A\,F_R(x,T;a)+B\,F_R(x,T;b)
\label{eq:two-step-finite}
\end{equation}
is built from the reflected distribution functions of Section~\ref{sec:boundary-reflection} for particles released at $a$ and at $b$. It solves the heat equation with the reflection rules the solver applies. Since each particle contributes its weight with probability $F_R(x,T;a)$ or $F_R(x,T;b)$, $m$ is also the mean of the reconstruction at the bin edges. The infinite-line reference is $u_L+A\,\Phi(x;a)+B\,\Phi(x;b)$, where $\Phi(x;\xi)$ is the Gaussian distribution function with mean $\xi$ and variance $2\alpha T$. This superposition of the error-function profiles~\eqref{eq:heat-exact} differs from $m$ by the domain truncation error. Comparing with each in turn separates the comparison-location error from the domain truncation error. Appendix~\ref{app:weighted} records the derivation, and a control that samples the final-time distribution directly, with no time stepping, reproduces $m$ to within the sampling scale of that control.

Let $U_N(b_j)$ denote one particle reconstruction at the right bin edge $b_j$, with bin width $h$, and let $g_j$ be the chosen deterministic reference value at the comparison point. For independent particles with fixed weights and allocation fractions $f_a=N_a/N$ and $f_b=N_b/N$, define $F_a(x)=F_R(x,T;a)$ and $F_b(x)=F_R(x,T;b)$, and write $m(x)=m(x,T)$ at the fixed final time. The expected squared error is
\begin{equation}
\begin{aligned}
\mathbb E[E_{\mathrm{total}}^2]
&=\mathbb E\big[\lVert U_N-g\rVert_{L_h^2}^2\big]=B_h^2+\frac{V_h}{N},\\
B_h^2&=h\sum_j\big(m(b_j)-g_j\big)^2,\\
V_h&=h\sum_j\left[\frac{A^2}{f_a}F_a(b_j)\big(1-F_a(b_j)\big)
+\frac{B^2}{f_b}F_b(b_j)\big(1-F_b(b_j)\big)\right].
\end{aligned}
\label{eq:weighted-identity}
\end{equation}
Here $B_h$ is the deterministic error floor, including the comparison-location and domain truncation effects and their cross term, and $V_h$ is the sampling coefficient. The expectation of the ensemble statistic $E_{\mathrm{total}}^2$ equals the expected squared error of one realization. Section~\ref{sec:analysis-cdf} explains the underlying empirical-distribution identity, and Appendix~\ref{app:weighted} derives the weighted form and the finite-ensemble factors used to predict bias and spread.

For $B_h>0$, the sampling and deterministic terms are equal at the transition count
\begin{equation}
N_\star=\frac{V_h}{B_h^2}.
\label{eq:transition-count}
\end{equation}
For a target root-mean-square error $\tau>B_h$, the expected squared error is at most $\tau^2$ when
\begin{equation}
N\ge\frac{V_h}{\tau^2-B_h^2}.
\label{eq:particle-count-rule}
\end{equation}
We round the count upward to an even integer to preserve the half-and-half allocation. When $B_h=0$, the rule reduces to $N\ge V_h/\tau^2$. For $\tau\le B_h$ and $V_h>0$, no finite particle count reaches the target.

Table~\ref{tab:two-step-predictions} gives the predictions. The predictions, accuracy target, and selected particle count were saved before the ensembles were run. The sampling coefficient $V_h$ varies only slightly across the bin counts, because it approximates a fixed integral of the two reflected distributions. The floors behave differently under bin refinement. Against the finite-interval reference at the bin edges the floor vanishes identically, since that reference is the mean of the reconstruction. Against the same reference at the bin centers it halves whenever $h$ halves, as expected from the first-order alignment error identified in Section~\ref{sec:heat-grid-paired}. Against the infinite-line reference at the bin edges it barely moves, because it measures the domain truncation error rather than a discretization. The transition count $N_\star$ of~\eqref{eq:transition-count} therefore quadruples with each halving of $h$ for the alignment floor and stays near $9\times10^3$ for the domain truncation floor.

\begin{table}[pos=htbp]
\centering\small
\caption{Predictions for the two-step problem from~\eqref{eq:weighted-identity}, computed before running any particles. $B_h$ is the predicted floor and $N_\star$ the transition count~\eqref{eq:transition-count} at which the sampling and deterministic contributions to the expected squared error are equal.}
\label{tab:two-step-predictions}
\begin{tabular}{@{}rccccccc@{}}
\toprule
& & & \multicolumn{2}{c}{finite interval, bin centers} & \multicolumn{2}{c}{infinite line, bin edges}\\
\cmidrule(lr){4-5}\cmidrule(lr){6-7}
$M$ & $h$ & $V_h$ & $B_h$ & $N_\star$ & $B_h$ & $N_\star$\\
\midrule
 50 & 0.0800 & 0.4488 & $2.255\times10^{-2}$ & 883    & $7.537\times10^{-3}$ & 7900\\
100 & 0.0400 & 0.4489 & $1.127\times10^{-2}$ & 3532   & $7.168\times10^{-3}$ & 8738\\
200 & 0.0200 & 0.4489 & $5.636\times10^{-3}$ & 14{,}130 & $6.991\times10^{-3}$ & 9185\\
400 & 0.0100 & 0.4489 & $2.818\times10^{-3}$ & 56{,}525 & $6.905\times10^{-3}$ & 9415\\
\bottomrule
\end{tabular}
\end{table}
\FloatBarrier

\begin{figure}[pos=htbp]
\centering
\includegraphics[width=1\textwidth]{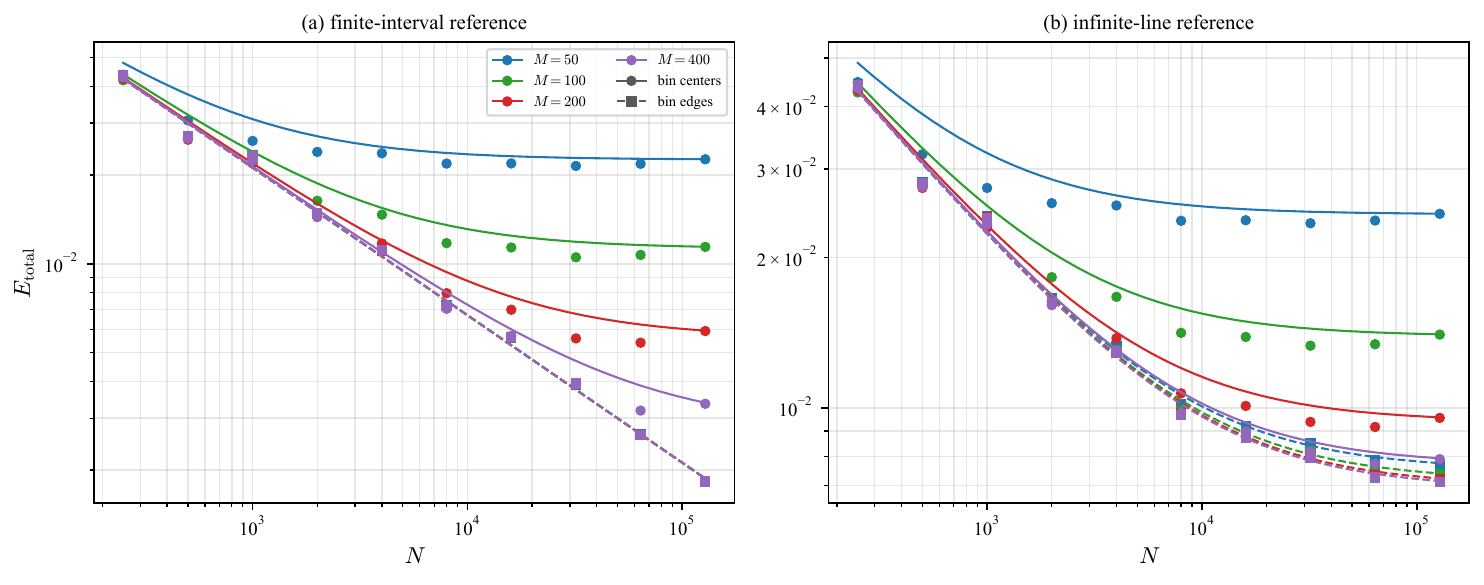}
\caption{Two-step heat problem. Measured $E_{\mathrm{total}}$ over thirty realizations (markers) against the predictions of~\eqref{eq:weighted-identity} (curves), for four bin counts, comparing at the bin centers (circles, solid) and at the bin edges (squares, dashed). (a)~against the finite-interval reference~\eqref{eq:two-step-finite}, where the bin-edge comparison has no floor. (b)~against the infinite-line reference, where the bin-edge floor is the domain truncation error and the bin-center floor also includes the alignment shift.}
\label{fig:two-step}
\end{figure}

Thirty realizations were then run at each of ten particle counts from $N=250$ to $N=128{,}000$, and each stored particle set was reconstructed at all four bin counts and compared with both references at both conventions. Figure~\ref{fig:two-step} shows the measured errors against the predicted curves. The predictions track the measurements over the whole range, through the sampling-dominated regime, the transition, and the floor. The four bin counts, two references, and two conventions are all evaluated on the same thirty trajectories at each particle count, so the $160$ resulting comparisons test the prediction across configurations rather than in independent experiments. Over those comparisons the ratio of measured to predicted $E_{\mathrm{total}}$ lies between $0.82$ and $1.10$, with median $0.98$, and the measured spreads lie between $0.88$ and $1.05$ of their predicted values.

\FloatBarrier

The deviations from prediction deserve a word, because they are correlated. Measured against the fluctuation scale that the covariance of Appendix~\ref{app:heat-uncertainty} predicts for thirty realizations, the bin-edge comparisons scatter about the prediction within $1.5$ standard deviations, while the bin-center comparisons fall below it at most particle counts, by up to three standard deviations. All ten particle counts reuse one list of thirty seeds, so their deviations are not independent. Repeating the experiment at $N=16{,}000$ on six independent blocks of thirty seeds gives measured-to-predicted ratios of $0.993$ at $M=50$ and $0.977$ at $M=200$, with standard deviations across the six blocks of $0.017$ and $0.040$, and the alignment term responsible for the offset changes sign between blocks. These blocks support sampling fluctuation as the explanation for the offset in the production ensembles, in the same way as the excess discussed in Appendix~\ref{app:heat-uncertainty}. The fluctuation scale itself comes from the Gaussian and delta-method approximations of that appendix.

Finally, we used the predictions to make an accuracy decision and tested it on seeds that had not been run. The target was a root-mean-square error of $8\times10^{-3}$ against the infinite-line reference at $M=200$. The model settles the choice of convention before any computation: comparing at the bin centers has floor $9.40\times10^{-3}$, so no particle count reaches the target, whereas comparing at the bin edges has floor $6.99\times10^{-3}$, and the rule~\eqref{eq:particle-count-rule} then requires $N\ge29{,}682$. Thirty fresh realizations at that count gave a measured error of $8.07\times10^{-3}$ against the predicted $8.00\times10^{-3}$, an agreement to within $1\%$ and within half a standard deviation of the prediction. At the same particle count the bin-center comparison gave $1.08\times10^{-2}$, above its own floor and above the target, as predicted. The rule constrains the expected squared error of a realization. The squared ensemble statistic fluctuates about that expectation, so a count chosen at the threshold leaves no margin and the measured error can fall on either side of the target. A user wanting the target met with margin would apply~\eqref{eq:particle-count-rule} with a tolerance below the target.

\section{Scalar Reaction--Diffusion Results}
\label{sec:fhn}

\subsection{Representative front propagation}
\label{sec:fhn-representative}

We solve the scalar reaction--diffusion equation with the FitzHugh--Nagumo-type cubic specified in Section~\ref{sec:fhn-benchmark}, using \(a=0.25\), \(D=0.5\), \(T=9\), \(\Delta t=0.01\), and \(N=500\)~globs on the domain~\([0,30]\), with the initial front location \(x_c=15\). Profiles are reconstructed at 500 uniform reconstruction points. This problem is the primary verification of the GRW method applied to a reaction--diffusion equation, where the weight update derived from the reaction term is active at each time step. These calculations use homogeneous Neumann boundary conditions, under which the weight is negated upon reflection. The front moves left from its initial location $x_c=15$ to $x\approx11.8$ by $T=9$, remaining well inside the domain $[0,30]$, so the boundary treatment has negligible influence on the results.
 
\begin{figure}[pos=htbp]
    \centering
\includegraphics[width=0.90\textwidth]{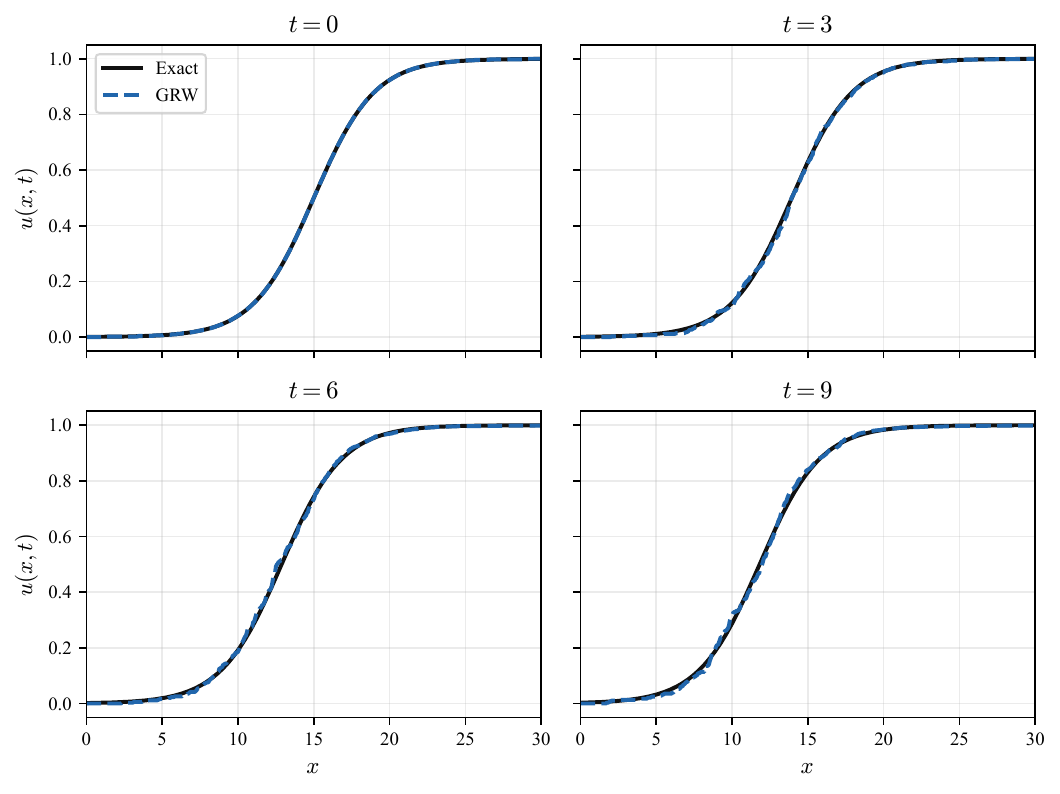}
    \caption{Scalar reaction--diffusion traveling front. GRW reconstruction versus exact traveling wave at \(t=0,3,6,9\), with \(a=0.25\), \(D=0.5\), \(N=500\). The GRW reconstruction (dashed line) and the exact front are so close that the two curves fall on top of each other over most of the domain.}
    \label{fig:fhn-comparison}
\end{figure}
\FloatBarrier
 
\begin{figure}[pos=htbp]
    \centering
    \includegraphics[width=1 \textwidth]{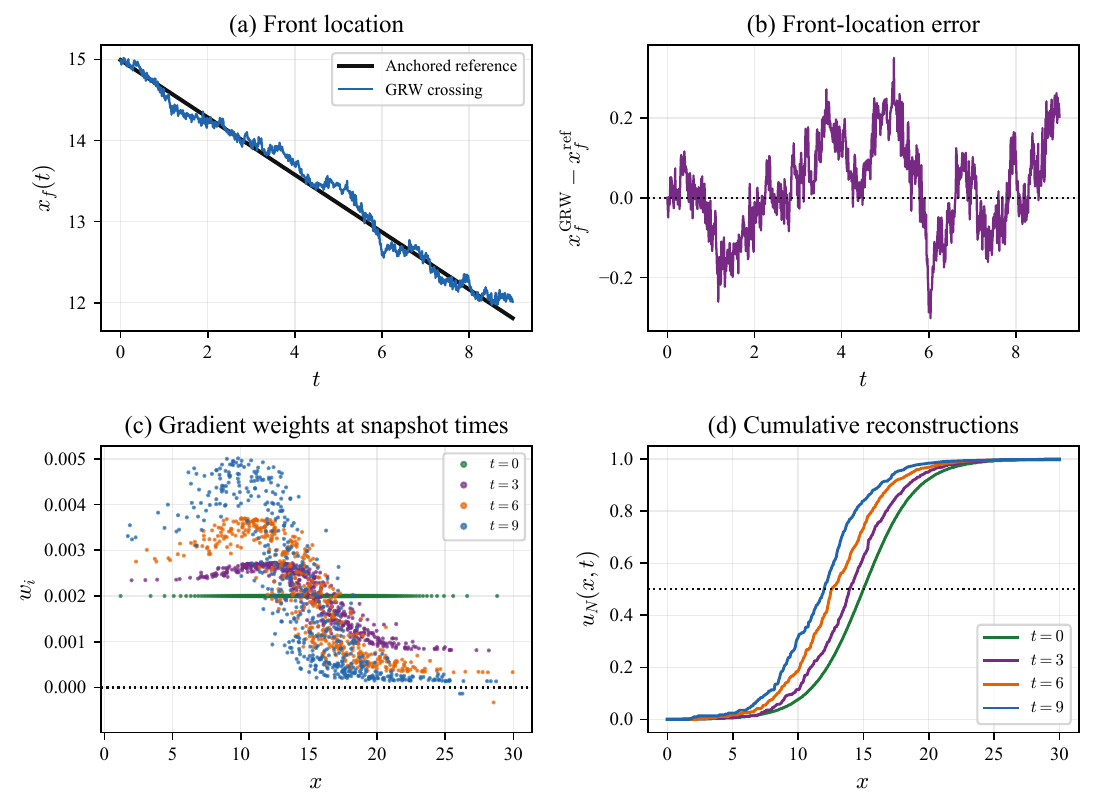}
    \caption{Scalar reaction--diffusion diagnostics. (a)~front location and (b)~front-location error versus time, (c)~signed glob weights and (d)~reconstructed fields at the snapshot times \(t=0,3,6,9\). In panels~(c) and~(d), colors distinguish the four snapshot times. \emph{Color is available in the online version of the article.}}
    \label{fig:fhn-diagnostics}
\end{figure}

Figure~\ref{fig:fhn-comparison} shows the front at $t=0,3,6,9$. The exact displacement over the run is $\theta T\approx3.18$. Figure~\ref{fig:fhn-diagnostics} resolves the front motion and underlying weight changes.

\FloatBarrier
 
The front-location comparison is anchored at the initial discrete crossing, which removes its initial offset as described in Section~\ref{sec:front-diagnostics}. At $T=9$, the location error is about $0.20$, or $6\%$ of the exact travel. A fit to the numerical trajectory gives speed $-0.345$, about $2.5\%$ below the exact magnitude $\theta\approx0.354$.

Initially every weight equals $1/N=0.002$. For $u\in[0,1]$, the tested reaction has $\max|R|=0.302$, so the multiplier at $\Delta t=0.01$ differs from one by at most $3.1\times10^{-3}$ and preserves weight sign. The two negative weights near $x\approx26$ at $t=9$ therefore arise from Neumann reflections, which negate the weights. The reconstructed profiles retain the sigmoid shape despite this redistribution. The following ensembles quantify profile and transport errors separately.

\subsection{Multi-seed profile, front location, speed, and aligned-profile convergence}
\label{sec:fhn-convergence}

The ensemble uses the same equation, domain $[0,30]$, homogeneous Neumann reflection, initial front location $x_c=15$, and time step as the representative problem, but shortens the final time from $T=9$ to $T=5$. The longer run demonstrates propagation over $3.18$ spatial units. The shorter run permits a thirty-seed refinement study while retaining $\theta T\approx1.77$ spatial units of travel. Thus the two studies test the same formulation for different purposes.

The convergence study uses
\begin{equation}
    D=0.5,
    \qquad
    a=0.25,
    \qquad
    x_c=15,
    \qquad
    T=5,
    \qquad
    \Delta t=0.01.
    \label{eq:fhn-study-parameters}
\end{equation}
With these values the wave speed parameter is
\begin{equation}
    \theta
    =
    \frac{1}{\sqrt{8}}
    \approx0.3536,
    \label{eq:fhn-theta}
\end{equation}
which gives the exact final front location
\begin{equation}
    x_{\mathrm{center}}^{\mathrm{exact}}(T)
    =
    15-\theta T
    \approx13.232.
    \label{eq:fhn-exact-center}
\end{equation}

Profile errors are measured against the exact logistic traveling wave~\eqref{eq:fhn-exact}, with every realization reconstructed at one fixed uniform set of $3001$ reconstruction points on $[0,L]$, where the norm~\eqref{eq:grid-l2} is evaluated. Because the imposed boundaries are homogeneous Neumann while the exact wave is an infinite-domain solution, the two differ slightly near the walls. We quantified this difference with a deterministic finite-difference solution of~\eqref{eq:fhn-pde} under the same Neumann conditions, converged in space and time step. That solution differs from the exact traveling wave by $8.3\times10^{-4}$ in $L_h^2$, which estimates the domain truncation error for this problem, the quantity that Section~\ref{sec:heat-grid-paired} computes exactly for the heat step. It is about a factor of $25$ below the smallest profile error in Table~\ref{tab:fhn-convergence}, so the exact wave serves as the reference. Center and speed errors need no numerical reference. They are measured against the exact front location~\eqref{eq:fhn-exact-center} and the exact velocity $-\theta$.

The numerical front location is extracted as in Section~\ref{sec:front-diagnostics}. Its average velocity is $c_{\mathrm{num}}=(x_{\mathrm{center}}^{\mathrm{num}}(T)-x_c)/T$. Since the exact location is $x_c-\theta T$,
\begin{equation}
|c_{\mathrm{num}}+\theta|=\frac{|x_{\mathrm{center}}^{\mathrm{num}}(T)-x_{\mathrm{center}}^{\mathrm{exact}}(T)|}{T}.
\label{eq:fhn-speed-center-relation}
\end{equation}
Location and speed errors differ by a fixed factor and therefore share the same log--log slope. We report the location-error confidence interval for both.

The fourth diagnostic, the aligned-profile error, separates front-shape error from translation error. The numerical profile is shifted horizontally by $x_{\mathrm{center}}^{\mathrm{ref}}-x_{\mathrm{center}}^{\mathrm{num}}$, the reference front location minus the extracted one, linearly interpolated onto the reference reconstruction points, and then compared with the reference in the discrete $L_h^2$ norm. Alignment removes the estimated translation, so the residual assesses shape error together with interpolation and finite-domain effects.

Thirty realizations were run for each particle count in
\begin{equation}
    N\in
    \left\{
        100,\,
        200,\,
        500,\,
        1000,\,
        2000,\,
        5000
    \right\},
    \label{eq:fhn-N-sequence}
\end{equation}
with the same sequence of seeds at each particle count. Table~\ref{tab:fhn-convergence} reports the mean direct-profile, front location, speed, and aligned-profile errors.

\begin{table}[pos=htbp]
    \centering
    \small
    \caption{Scalar reaction--diffusion GRW particle refinement ($D=0.5$, $a=0.25$, $T=5$, $\Delta t=0.01$, $S=30$). All four diagnostics decrease under refinement.}
    \label{tab:fhn-convergence}
    \begin{tabular}{rcccc}
        \toprule
        $N$
        & Profile $L_h^2$
        & Front-location error
        & Speed error
        & Aligned-profile $L_h^2$\\
        \midrule
        100  & 0.124 & 0.360 & 0.072 & 0.112\\
        200  & 0.081 & 0.195 & 0.039 & 0.070\\
        500  & 0.054 & 0.137 & 0.027 & 0.051\\
        1000 & 0.040 & 0.092 & 0.018 & 0.037\\
        2000 & 0.028 & 0.071 & 0.014 & 0.028\\
        5000 & 0.020 & 0.042 & 0.008 & 0.018\\
        \bottomrule
    \end{tabular}
\end{table}

\FloatBarrier

All four columns of Table~\ref{tab:fhn-convergence} decrease monotonically across the fiftyfold range of $N$. As noted in Section~\ref{sec:front-diagnostics}, the front-location diagnostic is local. In a single realization the crossing depends on the positions of a few nearby globs, so its seed-to-seed scatter is large. At $N=100$ the front-location error standard deviation is $0.26$ against a mean of $0.36$, and a monotone trend emerges only after averaging over realizations. A least-squares fit based on~\eqref{eq:loglog-rate-fit} gives the fitted exponents
\begin{subequations}
\label{eq:fhn-convergence-rates}
\begin{align}
    \text{profile }L_h^2
    &\sim N^{-0.460},
    &
    95\%\ \mathrm{CI}
    &=[-0.505,-0.415],
    \label{eq:fhn-profile-rate}\\
    \text{front-location error}
    &\sim N^{-0.524},
    &
    95\%\ \mathrm{CI}
    &=[-0.605,-0.445],
    \label{eq:fhn-center-rate}\\
    \text{speed error}
    &\sim N^{-0.524},
    &
    95\%\ \mathrm{CI}
    &=[-0.605,-0.445],
    \label{eq:fhn-speed-rate}\\
    \text{aligned-profile error}
    &\sim N^{-0.444},
    &
    95\%\ \mathrm{CI}
    &=[-0.479,-0.409].
    \label{eq:fhn-aligned-rate}
\end{align}
\end{subequations}
The front location and speed rows are identical by the exact relation~\eqref{eq:fhn-speed-center-relation}. The location, speed, and profile intervals contain the Monte Carlo rate $-1/2$, as do the heat intervals of Section~\ref{sec:heat-ensemble}. The aligned-profile interval lies above it, so that measure decays measurably more slowly than $-1/2$ over the tested range. All are fitted rates over a finite range of $N$.

The decrease in aligned-profile error shows that refinement improves the front shape after the estimated translation error is removed.

The measured total weight rises from $0.995$ at $N=100$ to $0.9997$ at $N=2000$, approaching the unit jump of the whole-line front. Section~\ref{sec:fhn-benchmark} distinguishes this whole-line balance from the finite-domain endpoint and reflection effects.

To give the rates a concrete scale, at $N=5000$ the speed error is about $0.008$, roughly $2.4\%$ of the exact speed magnitude $\theta$, and by~\eqref{eq:fhn-speed-center-relation} the front-location error is the same fraction of the distance $\theta T$ traveled over $T=5$.

The scalar reaction--diffusion update includes a reaction split, so unlike the heat case of Section~\ref{sec:brownian} its accuracy could in principle depend on the time step. At $N=2000$, varying $\Delta t$ from $0.04$ to $0.005$ gives nonmonotone profile errors between $0.0265$ and $0.0325$, with no resolved time-step trend relative to the seed-to-seed variability. At these parameters the time step is therefore not the limiting factor.

\begin{figure}[pos=htbp]
    \centering
    \includegraphics[width=1\textwidth]{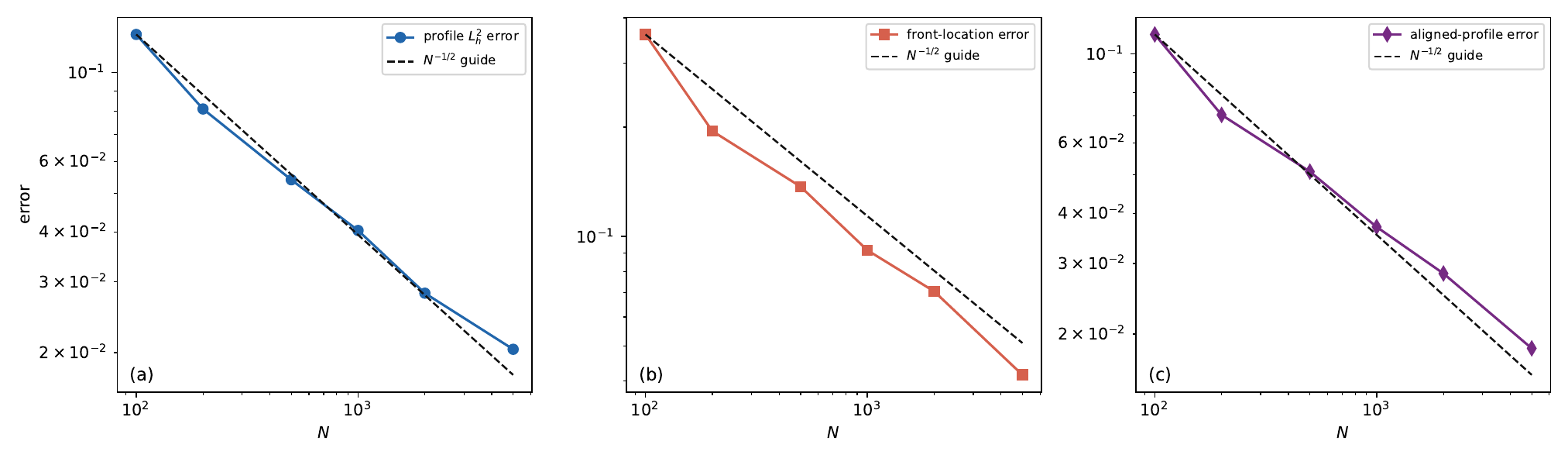}
    \caption{Scalar reaction--diffusion convergence under particle refinement over the thirty-seed ensemble. (a)~mean profile $L_h^2$ error, (b)~front-location error, (c)~aligned-profile $L_h^2$ error, each with the $N^{-1/2}$ guide. Fitted slopes and realization-level confidence intervals are given in~\eqref{eq:fhn-convergence-rates}.}
    \label{fig:fhn-convergence}
\end{figure}

Figure~\ref{fig:fhn-convergence} plots the profile, location, and aligned-profile errors against the $N^{-1/2}$ guide. Speed error is omitted because it is location error divided by $T$.

\FloatBarrier

\section{Cole--Hopf Burgers Results and Diagnosis}
\label{sec:colehopf}

\subsection{Representative recovery and error decomposition}
\label{sec:burgers-representative}

We now turn to the Burgers problem of Section~\ref{sec:burgers-benchmark}, which uses \(A=1\), \(\nu=0.5\), \(T=0.5\), \(\Delta t=0.005\), and \(P=M=400\) on the domain~\([0,4]\), with the shock centered at \(x_c=L/2=2\). The \(P\) initialization points produce \(P-1=399\) \(\phi_x\)-globs, and the bin width \(h\approx0.01\) is the same here and in the domain study below (where \(P=M=100L\)). In this representative run, the reconstructed~\(\phi_N\) reaches a minimum near \(0.343\), matching the closed-form value \(0.341\) of Section~\ref{sec:burgers-benchmark} and staying well above the denominator floor of \(0.133\), so the floor is never active.

\begin{figure}[pos=htbp]
    \centering
\includegraphics[width=0.75\textwidth]{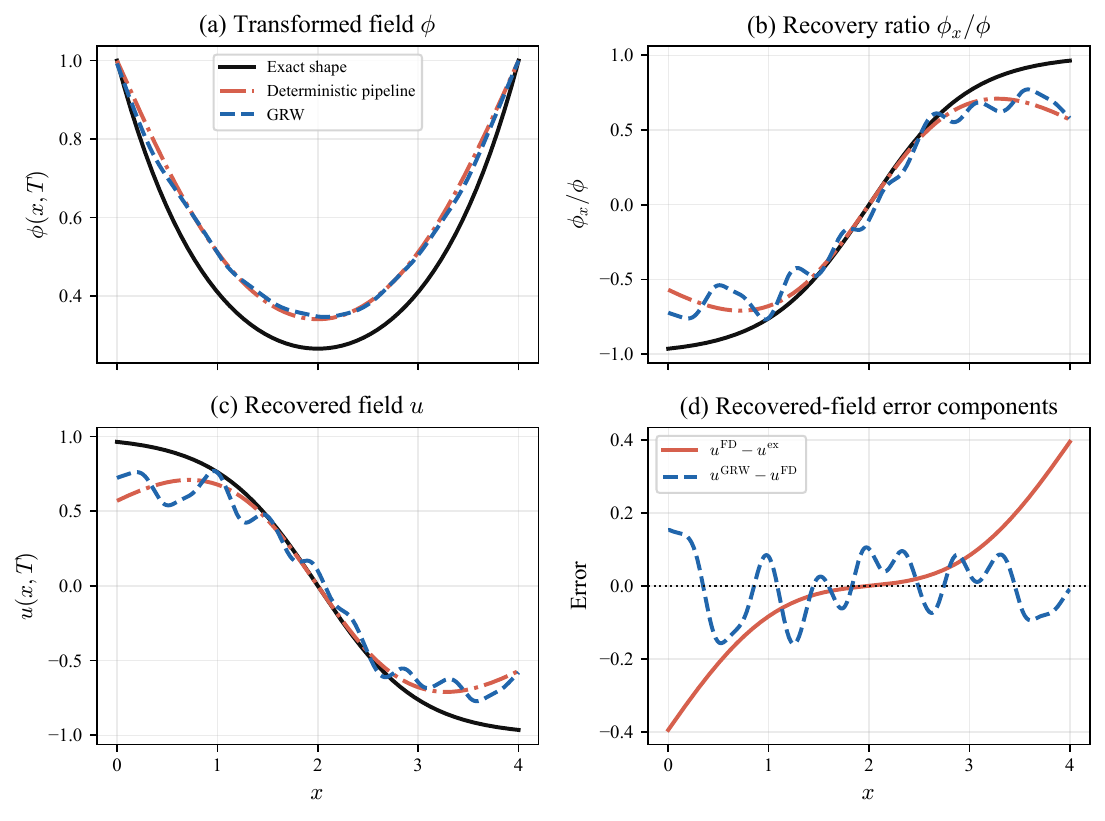}
    \caption{Cole--Hopf Burgers diagnostics at \(T=0.5\), \(P=M=400\), \(L=4\). (a)~transformed field~\(\phi\), (b)~recovery ratio~\(\phi_x/\phi\), (c)~recovered field~\(u\), (d)~pointwise error components \(u^{\mathrm{FD}}-u^{\mathrm{ex}}\) and \(u^{\mathrm{GRW}}-u^{\mathrm{FD}}\), whose RMSEs are \(E_{\mathrm{det}}\) and \(E_{\mathrm{GRW}}\). The recovery error concentrates at the domain endpoints, with the deterministic component \(E_{\mathrm{det}}\) the larger of the two at this domain size.}
    \label{fig:burgers-diagnostics}
\end{figure}
\FloatBarrier
 
Figure~\ref{fig:burgers-diagnostics} shows how error enters at each stage of the Cole--Hopf procedure. Panel~(a) compares the GRW solution for~\(\phi\) with both the deterministic finite-domain solution of the transformed heat equation~\eqref{eq:cole-hopf-heat} and the exact infinite-domain profile. The deterministic solution already departs from the exact profile near the boundaries, before any particle noise is introduced. Panel~(b) shows the recovery ratio~\(\phi_x/\phi\). Small differences in~\(\phi\) are amplified by the differentiation and division in the recovery formula~\eqref{eq:cole-hopf-recovery}. The division is most sensitive near the shock center, where~\(\phi\) is smallest. The largest errors in the recovered field, however, occur at the domain endpoints, where the finite-domain treatment holds the boundary values of~\(\phi\) fixed and the derivative~\(\phi_x\) must be formed from one-sided differences. Panels~(c) and~(d) show the recovered fields and their error components. Both recovered fields depart from the exact shock near the endpoints, and $E_{\mathrm{det}}$ is the larger component at $L=4$.

\subsection{Domain-size study}
\label{sec:results-burgers}

To examine how the error components depend on domain size, we vary $L=4,6,8,10$ with $P=M=100L$, so that the particle density is fixed at $100$ points per unit length while the bin width stays near $h=0.010$ and the smoothing width varies only from $\sigma_x=0.1203$ to $0.1201$ as the domain widens. The deterministic component $E_{\mathrm{det}}$ is a property of the finite-domain solve and is computed once per domain. The particle and total errors are computed over a thirty-seed ensemble at each domain, following Section~\ref{sec:ensemble-vocab}, and Table~\ref{tab:burgers-domain} reports their means and seed-to-seed standard deviations.

\begin{table}[pos=htbp]
\centering
\small
\caption{Burgers error decomposition versus domain size with $P=M=100L$, evaluated at $T=0.5$ over $S=30$ realizations per domain. Entries are the RMSE quantities of~\eqref{eq:burgers-error-decomposition}. The deterministic component carries no seed dependence. The last column restates $E_{\mathrm{det}}$ in the $L_h^2$ norm, in which it is essentially independent of $L$.}
\label{tab:burgers-domain}
\begin{tabular}{@{}cccccc@{}}
\toprule
$L$ & $P=M$ & $E_{\mathrm{det}}$ & $E_{\mathrm{GRW}}$ & $E_{\mathrm{total}}$ & $E_{\mathrm{det}}$ ($L_h^2$)\\
\midrule
4  & 400  & 0.172 & $0.104\pm0.021$ & $0.200\pm0.021$ & 0.344\\
6  & 600  & 0.138 & $0.127\pm0.023$ & $0.190\pm0.024$ & 0.339\\
8  & 800  & 0.120 & $0.140\pm0.020$ & $0.187\pm0.025$ & 0.339\\
10 & 1000 & 0.107 & $0.160\pm0.025$ & $0.193\pm0.022$ & 0.338\\
\bottomrule
\end{tabular}
\end{table}
\FloatBarrier

The deterministic component decreases in RMSE from $0.172$ at $L=4$ to $0.107$ at $L=10$. By the norm conversion of Section~\ref{sec:norms}, the deterministic error in the $L_h^2$ norm is $E_{\mathrm{det}}\sqrt{hM}$, which gives $0.344$, $0.339$, $0.339$, and $0.338$ across the four domains, essentially unchanged by domain size. The RMSE decrease is therefore the $\sqrt{hM}$ factor of that conversion, not a reduction in the error itself. The boundary mismatch stays concentrated near the endpoints (panel~(d) of Figure~\ref{fig:burgers-diagnostics}), so widening the domain dilutes it in a whole-domain average without reducing it. The particle component increases from $0.104\pm0.021$ to $0.160\pm0.025$, consistent with the larger dynamic range of the normalized transformed field on wider domains. Their ordering reverses between $L=6$ and $L=8$, as Table~\ref{tab:burgers-domain} shows. The mean total error remains between $0.187$ and $0.200$. A linear fit versus $L$ (see Section~\ref{sec:ensemble-bootstrap}) has slope $-1.20\times10^{-3}$ with $95\%$ interval $[-3.02\times10^{-3},7.30\times10^{-4}]$, which contains zero. Widening the domain therefore shifts the error composition, while the total-error trend remains unresolved over the tested sizes.

\subsection{Controlled attribution of the accuracy floor}
\label{sec:colehopf-diagnosis}

The Cole--Hopf construction has two discretizations, the $P$ initialization points and the $M$ reconstruction bins of Section~\ref{sec:grid-free}. The coupled design of Section~\ref{sec:burgers-benchmark} ties them together through $P=M$, and the controls below separate them. Under coupled refinement, the mean of the realization errors in the $L_h^2$ norm~\eqref{eq:grid-l2} is
\begin{center}\small
\begin{tabular}{@{}lcccccc@{}}
\toprule
$P=M$ & 50 & 100 & 200 & 400 & 800 & 1600\\
\midrule
error & 1.650 & 0.508 & 0.426 & 0.391 & 0.421 & 0.408\\
\bottomrule
\end{tabular}
\end{center}
The error falls steeply through $P=M=100$ and then holds near $0.41$.

The plateau can be compared with the domain study total error through the norm conversion of Section~\ref{sec:norms}. At $L=4$,
\[
\lVert e\rVert_{L_h^2}=\sqrt{hM}\,\operatorname{RMSE}(e)\approx 2\,\operatorname{RMSE}(e).
\]
The ensemble-mean total error of Table~\ref{tab:burgers-domain} at $L=4$ is then $0.200$ in RMSE and $0.401$ in the $L_h^2$ norm~\eqref{eq:grid-l2}. The plateau is thus the ensemble total error of the representative configuration, and the controls below attribute it.

The coupled design varies three numerical choices at once, the number of $\phi_x$-globs, the number of reconstruction bins, and the physical smoothing bandwidth. Because the kernel width is fixed at $12$ bins, the physical standard deviation is
\[
\sigma_x=\frac{12L}{M-1}.
\]
This shrinks from about $0.98$ at $P=M=50$, comparable to the shock transition width $2\nu/A=1$, to $0.03$ at $P=M=1600$. Each control below holds two of these choices fixed, or removes the particles entirely.

\subsubsection{Boundary control}
\label{sec:colehopf-boundary}
The boundary control tests how much of the deterministic component in Table~\ref{tab:burgers-domain} comes from the transformed endpoint data. The transformed solve of Section~\ref{sec:burgers-decomposition} holds the endpoint values of $\phi$ fixed at their initial values, whereas the exact transformed field $\phi(x,t)=e^{\nu k^{2}t}\phi_0(x)$ derived in Section~\ref{sec:burgers-benchmark} grows its endpoint values by the factor $e^{\nu k^{2}T}\approx1.28$ over the run.

\begin{figure}[pos=htbp]
    \centering
\includegraphics[width=1.00\textwidth]{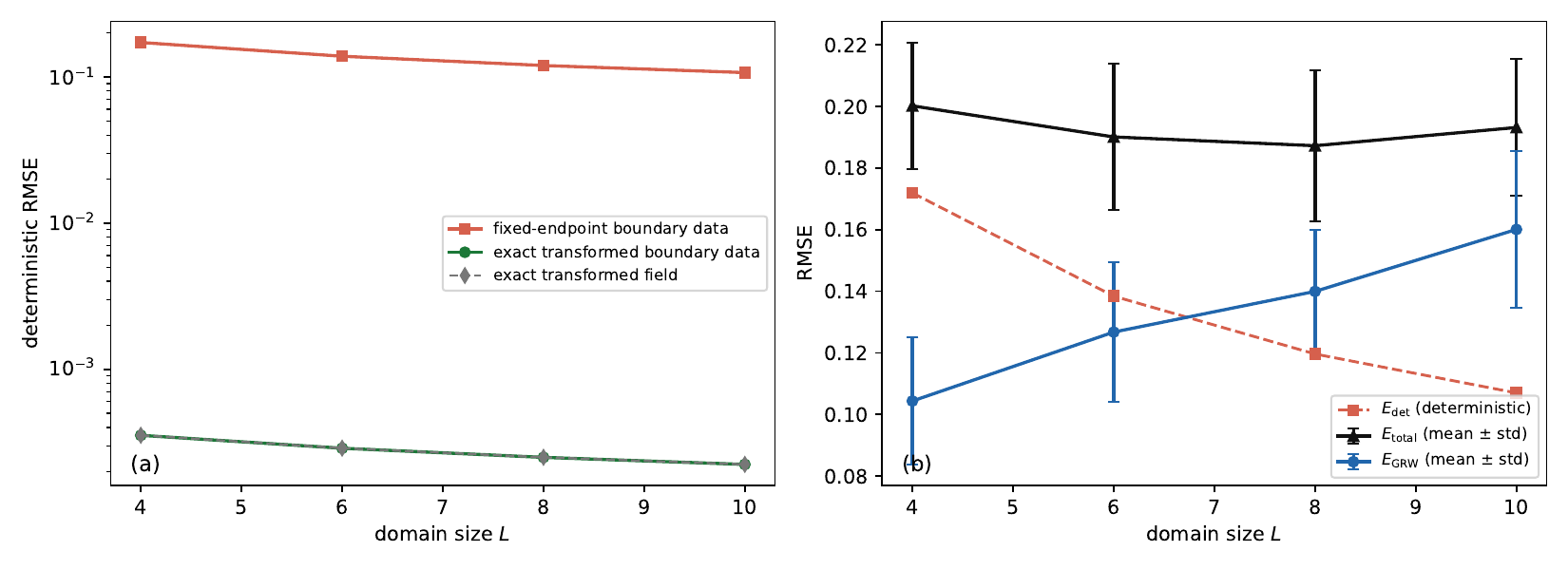}
    \caption{(a)~Deterministic boundary controls versus domain size. The fixed-endpoint solve carries the dominant deterministic error, while the solve with exact transformed boundary data falls to the level of the pure inversion control. (b)~The thirty-seed domain decomposition of Table~\ref{tab:burgers-domain}, with opposing component trends and a total that shows no monotone trend across the tested domains. Both panels use RMSE. Table~\ref{tab:burgers-domain} also gives the deterministic $L_h^2$ error, which is essentially unchanged as the domain widens.}
    \label{fig:burgers-boundary-domain}
\end{figure}

Supplying the same finite-difference solve with these exact time-dependent values (verified against the closed form to $8\times10^{-7}$) reduces the deterministic error at both domains:
\[
0.172\longrightarrow 3.5\times10^{-4}\ (L=4),\qquad 0.107\longrightarrow 2.2\times10^{-4}\ (L=10),
\]
as shown in Figure~\ref{fig:burgers-boundary-domain}(a). For this stationary shock, the exact endpoints remove all but $0.2\%$ of the deterministic error at $L=4$, a factor of about $490$.

The particle calculations retain the original fixed endpoint values of $\phi$. Section~\ref{sec:analysis-colehopf} relates this boundary control to the physical Burgers boundary data.

\subsubsection{Inversion control}
\label{sec:colehopf-inversion}
Supplying the accurate transformed field directly to the same differentiation and inversion procedure yields
\[
7\times10^{-4}\ \text{in } L_h^2\qquad (3.5\times10^{-4}\ \text{RMSE}),
\]
the same level reached by the boundary-consistent solve. The discrete recovery error is negligible at the scale of the boundary effect in this control.

\FloatBarrier

\subsubsection{Decoupled refinement}
\label{sec:colehopf-decoupled}
\begin{figure}[pos=htbp]
    \centering
\includegraphics[width=1.00\textwidth]{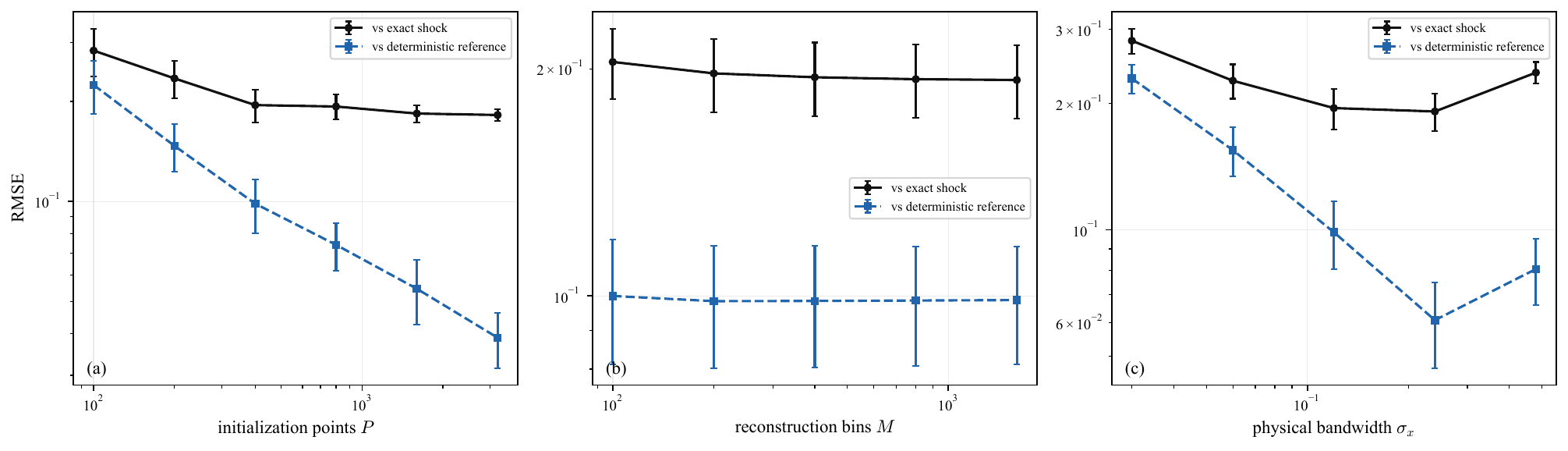}
    \caption{Decoupled refinement of the Cole--Hopf GRW procedure over twenty-seed ensembles (RMSE mean $\pm$ one standard deviation), against the exact shock and the deterministic fixed-endpoint reference. (a)~Initialization-point refinement at fixed reconstruction bins and bandwidth. (b)~Reconstruction-bin refinement at fixed initialization points and bandwidth. (c)~Bandwidth sweep at fixed initialization points and reconstruction bins.}
    \label{fig:burgers-decoupled}
\end{figure}

Table~\ref{tab:burgers-decoupled-design} lists the three decoupled designs. Each holds two of the three numerical choices fixed and varies the third over a twenty-seed ensemble. Figure~\ref{fig:burgers-decoupled} reports the RMSE against both the exact shock and the deterministic fixed-endpoint reference, with the latter comparison measuring particle sampling, initialization, and reconstruction effects.

\begin{table}[pos=htbp]
\centering
\small
\caption{Design of the decoupled refinement study. $P$ is the number of initialization points ($P-1$ globs), $M$ the number of reconstruction bins, and $\sigma_x$ the physical smoothing bandwidth. Every configuration uses the same twenty seeds.}
\label{tab:burgers-decoupled-design}
\begin{tabular}{@{}lccc@{}}
\toprule
Control & $P$ & $M$ & $\sigma_x$\\
\midrule
initialization-point refinement & $100$--$3200$ & $400$ & $0.12$\\
reconstruction-bin refinement & $400$ & $100$--$1600$ & $0.12$\\
bandwidth sweep & $400$ & $400$ & $0.03$--$0.48$\\
\bottomrule
\end{tabular}
\end{table}
\FloatBarrier

At fixed $M=400$ and $\sigma_x=0.12$, refining the initialization points from $P=100$ to $3200$ reduces the particle RMSE from $0.224$ to $0.039$. The fitted exponent is $-0.495$, with $95\%$ confidence interval $[-0.521,-0.470]$. Because $P$ also sets the trapezoidal quadrature of the transformed initial condition~\eqref{eq:cole-hopf-initialization}, this sweep refines the initialization points and the particle count together, and its discrepancy includes binning, smoothing, and particle noise. Its reference is the fixed-endpoint deterministic solve, so this finding is independent of the boundary control in Section~\ref{sec:colehopf-boundary}. At fixed $P=400$ and bandwidth, refining the reconstruction bins from $M=100$ to $1600$ changes the particle RMSE from $0.100$ to $0.0987$, about $1.3\%$, showing weak reconstruction-bin dependence over this range.

The bandwidth sweep at $P=M=400$ has its lowest tested particle RMSE, $0.061$, at $\sigma_x=0.24$. The sampling spread decreases from $0.224$ at $\sigma_x=0.03$ to $0.032$ at $\sigma_x=0.48$ across the sweep, while the ensemble-mean discrepancy from the deterministic reference rises to $0.076$ at $\sigma_x=0.48$. This measured bias--variance tradeoff helps explain why coupling the physical bandwidth to particle refinement can obscure the particle-count dependence.

\FloatBarrier

\subsubsection{Amplification control}
\label{sec:colehopf-amplification}
\begin{figure}[pos=htbp]
    \centering
\includegraphics[width=0.90\textwidth]{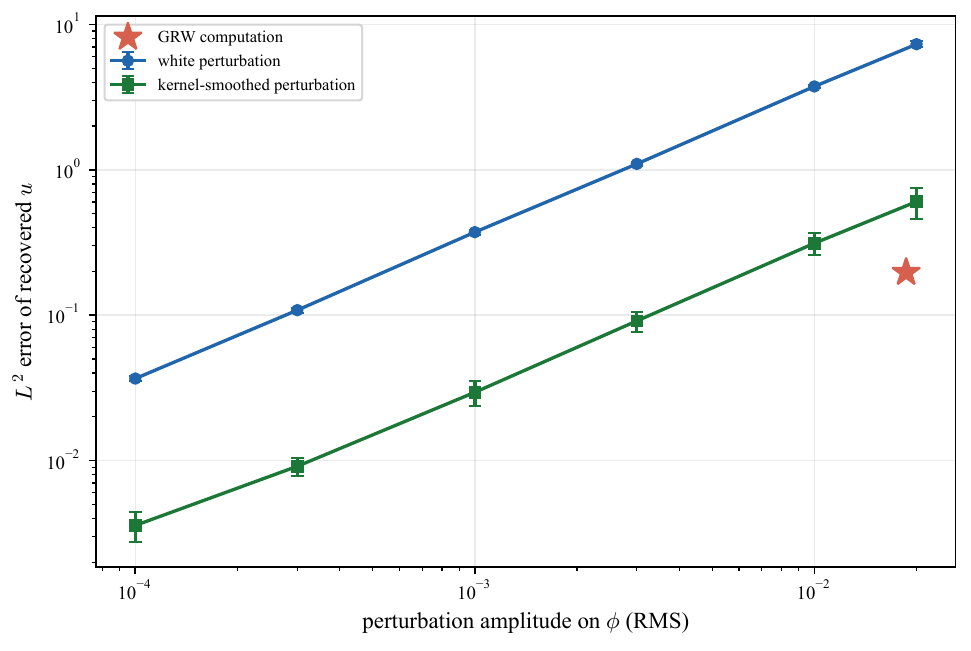}
    \caption{Response of the recovery map to perturbations of the transformed field, twenty realizations per amplitude, for white and kernel-smoothed perturbations of equal root-mean-square amplitude. Recovered-field errors are reported in $L_h^2$. The star marks the measured Cole--Hopf GRW calculation.}
    \label{fig:burgers-perturbation-response}
\end{figure}

We apply additive perturbations of chosen root-mean-square amplitude to the exact transformed profile normalized to unit endpoint values, and recover $u$ by numerical differentiation and division, using twenty realizations per amplitude. This control clips the perturbed $\phi$ only at $10^{-10}$, while the particle solver uses the initial-profile floor of Section~\ref{sec:burgers-benchmark}. The white family perturbs the transformed field independently at each reconstruction bin. The smoothed family passes the same perturbations through the Gaussian kernel of the reconstruction and rescales them to the same amplitude, so the two families differ only in their spatial correlation. Because the particle error is itself kernel-smoothed, the smoothed family models its structure, and the white family shows the response to uncorrelated error. Figure~\ref{fig:burgers-perturbation-response} shows the recovered-field response. At amplitudes $10^{-3}$, $10^{-2}$, and $2\times10^{-2}$, the two families produce
\[
\text{white } (0.373,\ 3.75,\ 7.31),\qquad \text{smoothed } (0.030,\ 0.312,\ 0.604).
\]
The white response is about $12$ times the smoothed response at each of the three amplitudes shown, so the response depends strongly on the perturbation's spatial structure. The measured GRW reconstruction (transformed-field RMSE $1.9\times10^{-2}$ and particle error $0.197$ in $L_h^2$, comparable to the $E_{\mathrm{GRW}}$ of Table~\ref{tab:burgers-domain} after the norm conversion) lies below the kernel-smoothed response and far below the white-noise response. Section~\ref{sec:analysis-colehopf} explains this dependence through the derivative of the relative error in the transformed field.

\FloatBarrier

Together, the controls explain the plateau through the transformed boundary data and the reconstruction choices varied in coupled refinement. They also identify two separate improvements to pursue, consistent transformed boundary data and selection of a physical smoothing bandwidth. The identities in Section~\ref{sec:analysis-colehopf} explain how each affects the recovered Burgers field.

\Needspace{6\baselineskip}
\section{\texorpdfstring{Analytical basis and scope of the attributions}{Analytical basis and scope of the attributions}}
\label{sec:analysis}
We now state the analytical results underlying the controls of Sections~\ref{sec:heat}--\ref{sec:colehopf}. They concern cumulative reconstruction, empirical-distribution error, reaction and boundary compatibility, and Cole--Hopf sensitivity. Their assumptions determine when the diagnostics apply to other problems.

\subsection{What a cumulative reconstruction represents}
\label{sec:analysis-coordinates}
Let the reconstruction bins have edges $b_0<b_1<\dots<b_M$ and width $h$, let $c_j=(b_{j-1}+b_j)/2$ be the bin centers, and let $m_j$ be the total signed weight of the globs in $(b_{j-1},b_j]$. The bin-and-sum reconstruction is
\begin{equation}
U_j=u_L+\sum_{\ell=1}^{j}m_\ell ,
\label{eq:cumulative-bins}
\end{equation}
which counts every glob to the left of $b_j$ and is therefore the reconstruction of $u$ at the right edge $b_j$. This holds for signed weights as well as for the equal positive weights of the heat step. The paired study of Section~\ref{sec:heat-grid-paired} compares $U_j$ with $u^{\mathrm{ref}}(b_j)$ and with $u^{\mathrm{ref}}(c_j)$. If $u^{\mathrm{ref}}$ has a bounded second derivative, a Taylor expansion about $c_j$ with $b_j-c_j=h/2$ gives
\begin{equation}
u^{\mathrm{ref}}(b_j)-u^{\mathrm{ref}}(c_j)=\frac{h}{2}\,u^{\mathrm{ref}}_x(c_j)+O(h^2).
\label{eq:halfbin}
\end{equation}
Comparing $U_j$ with the value at the center therefore misplaces the comparison by half a bin and adds an error of first order in $h$, proportional to $u_x$. This shift is the leading contribution to the bin-center floor in Table~\ref{tab:heat-grid-paired} and is concentrated where the profile is steep. The domain truncation error and its cross term with the shift explain the departure of the floor ratio from the bin-width ratio in Section~\ref{sec:heat-grid-paired}. Identity~\eqref{eq:paired-identity} relates two such comparisons of the same reconstruction. The same edge interpretation applies to any method that recovers a field by summing binned signed masses, and reconstructions of other types must be interpreted on their own terms. The coupled treatment of Section~\ref{sec:heat-ensemble} evaluates the cumulative sum directly at its own points, so no shift arises, and a density obtained by dividing $m_j$ by $h$ is a bin average rather than an edge value.

\subsection{Expected error of the empirical distribution function}
\label{sec:analysis-cdf}
For the heat step the globs carry equal weight $1/N$, do not interact, and move independently. With $u_L=0$ and a unit jump, the reconstruction~\eqref{eq:cumulative-bins} is then the empirical distribution function
\[
F_N(x)=\frac1N\sum_{i=1}^N I_i(x),\qquad I_i(x)=\mathbf 1_{\{X_i\le x\}},
\]
of $N$ independent particle positions, each with the same distribution function $F$. Here $F=F_R$ is the distribution function of a reflected Brownian particle on $[0,L]$, introduced in Section~\ref{sec:boundary-reflection} and given explicitly in Appendix~\ref{app:heat-uncertainty}. The folded Gaussian steps of Section~\ref{sec:boundary-reflection} sample this distribution exactly at every step. Symmetry of the Gaussian increments makes repeated folding equivalent in distribution to folding the unreflected Gaussian walk once at the final time. Appendix~\ref{app:heat-uncertainty} proves this equivalence and gives the resulting Neumann heat kernel~\citep{carslawjaeger1959,pavliotis2014}. The time step therefore contributes no error to the particle distribution in this problem. This statement is specific to constant diffusivity and simple reflection, and the reaction split of Section~\ref{sec:fhn} introduces a separate time-step dependence.

Since each indicator $I_i(x)$ is a Bernoulli variable with success probability $F(x)$, and the particles are independent, $\mathbb E[F_N(x)]=F(x)$ and $\operatorname{Var}(F_N(x))=F(x)(1-F(x))/N$~\citep{vandervaart1998}. Bertaglia et al.~\citep{bertaglia2024} use the same property for gradient-based Monte Carlo, giving the exact pointwise second moment of the error of an empirical distribution function on the whole line together with a bound on its $L^p$ norm. We apply this variance in the reported discrete norm, using the reflected finite-interval distribution and recording the comparison-location error in the deterministic term. Applying $\mathbb E[(Y-g)^2]=\operatorname{Var}(Y)+(\mathbb EY-g)^2$ at each comparison point and summing gives, for deterministic comparison values $g_j$,
\begin{equation}
\mathbb{E}\Big[h\sum_{j}\big(F_N(b_j)-g_j\big)^2\Big]
=h\sum_j\Big[\frac{F(b_j)\big(1-F(b_j)\big)}{N}+\big(F(b_j)-g_j\big)^2\Big].
\label{eq:cdf-identity}
\end{equation}
At fixed comparison points, the sampling term scales exactly as $N^{-1}$. The deterministic term is the squared systematic error. Its square root gives the bin-center floor of Section~\ref{sec:heat-grid-paired} when $g_j=u^{\mathrm{ex}}(c_j)$, and the domain truncation error when $g_j=u^{\mathrm{ex}}(b_j)$. Correlations between different comparison points do not enter~\eqref{eq:cdf-identity}, because the left side is a sum of pointwise squares, although they do enter the variance of that sum (Appendix~\ref{app:heat-uncertainty}). The identity gives the expected squared error of a single realization. Multiplying its sampling term by $(S-1)/S$ gives $\mathbb E[E_{\mathrm{spread}}^2]$ for an ensemble of $S$ realizations, as derived in Appendix~\ref{app:heat-uncertainty}. Appendix~\ref{app:heat-uncertainty} compares the measured spreads with these root-mean-square predictions and assesses their fluctuations using a Gaussian approximation and independent seed blocks.

For a step of amplitude $u_R-u_L$ the reconstruction is $u_L+(u_R-u_L)F_N$, and, for a physical reference $u_L+(u_R-u_L)g_j$, both sides of~\eqref{eq:cdf-identity} scale by $(u_R-u_L)^2$.

For fixed signed or unequal weights, the same pointwise argument adds the independent variance contributions with squared weights. Section~\ref{sec:heat-two-step} states the resulting identity~\eqref{eq:weighted-identity} for the two-step study, and Appendix~\ref{app:weighted} derives its finite-ensemble factors. In the reaction problem, the weights depend on the evolving reconstructed field and interact. For Cole--Hopf, the weighted argument describes the unsmoothed cumulative reconstruction of $\phi$, while predicting error in the recovered $u$ also requires the smoothing, weight correction, and nonlinear recovery.

\subsection{Reaction laws and boundary data for the weight update}
\label{sec:analysis-reaction}
The weight update of Section~\ref{sec:reaction-update} applies to the class of scalar reaction--diffusion equations
\[
u_t=Du_{xx}+f(u),\qquad D>0 \text{ constant},
\]
with $f$ continuously differentiable and $u$ smooth enough for the equation to be differentiated. Then $v=u_x$ satisfies $v_t=Dv_{xx}+f'(u)v$, and the multiplier $1+\Delta t\,f'(u_i^{*})$ in~\eqref{eq:reaction-weight-update} is an explicit reaction update. The probability restriction for its single-event birth--death interpretation is stated in Section~\ref{sec:reaction-update}. Stability and convergence require separate analysis. Table~\ref{tab:reaction-family} lists the derivatives required by the update, including the classical Fisher--KPP and bistable Nagumo particle problems studied by Puckett~\citep{puckett1989} and Mascagni~\citep{mascagni1995}. The last column identifies the cases computed here.

\begin{table}[pos=htbp]

\centering\small
\caption{Reaction laws covered by the weight update~\eqref{eq:reaction-weight-update}. The family parameter $\rho$ is positive. The last column gives the section in which the law is computed in this paper.}
\label{tab:reaction-family}
\begin{tabular}{@{}llll@{}}
\toprule
Equation & $f(u)$ & $f'(u)$ & Computed in\\
\midrule
Heat & $0$ & $0$ & Section~\ref{sec:heat}\\
Fisher--KPP & $\rho u(1-u)$ & $\rho(1-2u)$ & ---\\
Bistable Nagumo, $0<a<1$ & $\rho u(1-u)(u-a)$ & $\rho\big[-3u^2+2(1+a)u-a\big]$ & ---\\
Monostable cubic~\eqref{eq:fhn-reaction} & $\tfrac{D}{2}u(1-u)(u-u_\star)$, $u_\star<0$ & Eq.~\eqref{eq:fhn-statistic} & Section~\ref{sec:fhn}\\
\bottomrule
\end{tabular}
\end{table}

Boundary data are the other condition. At an endpoint $x=a$ with prescribed value $u(a,t)=g(t)$, evaluating the equation at the wall gives
\begin{equation}
D\,v_x(a,t)=g'(t)-f\big(g(t)\big)
\label{eq:boundary-compat}
\end{equation}
whenever $u$ is smooth up to the boundary. The weight-preserving reflection of Section~\ref{sec:boundary-reflection} extends $v$ evenly across the wall and so imposes $v_x(a,t)=0$. It is compatible with the prescribed value only when $g'(t)=f(g(t))$, and for a constant boundary state only when $f(g)=0$, i.e., when the boundary value is an equilibrium of the reaction. The heat step and the original fixed-endpoint calculation for the transformed field $\phi$ use constant boundary states and $f=0$. The reaction front uses weight-negating reflection, which extends $v$ oddly and imposes $v(a,t)=0$, the homogeneous Neumann condition, for which no compatibility condition on $g$ arises. Other time-dependent Dirichlet data require an additional boundary treatment. Equation~\eqref{eq:boundary-compat} states continuum compatibility for these reflection rules, and the accuracy of their finite-particle enforcement is a separate discretization question. The update also has clear limits. A spatially varying reaction $f(x,u)$ adds a term $f_x$ to the equation for $v$ that is not a multiplier of $v$. For the nondivergence equation $u_t=D(x)u_{xx}+f(u)$, differentiation gives $v_t=Dv_{xx}+D_xv_x+f'(u)v$, whereas the divergence-form equation $u_t=(D(x)u_x)_x+f(u)$ gives $v_t=(Dv)_{xx}+f'(u)v$. Systems couple the weights of several fields, and in more than one dimension $v$ is a vector whose recovery requires a potential or a line integral. Each of these needs structure beyond the update used here.

\subsection{Sensitivity of Cole--Hopf recovery}
\label{sec:analysis-colehopf}
The recovery $u=-2\nu\,\partial_x\log\phi$ acts on the logarithm of the transformed field, so its error depends on the relative error of $\phi$ and on the derivative of that relative error. For continuously differentiable $\phi>0$ and $\varepsilon$ on the finite interval, let $\widetilde\phi=\phi(1+\varepsilon)$ with $1+\varepsilon>0$. Subtracting the two logarithmic derivatives gives
\begin{equation}
\widetilde u-u=-2\nu\,\partial_x\log(1+\varepsilon)=-2\nu\,\frac{\varepsilon_x}{1+\varepsilon}
\label{eq:colehopf-perturbation}
\end{equation}
exactly. If $\lVert\varepsilon\rVert_\infty<1$, then $|1+\varepsilon(x)|\ge1-\lVert\varepsilon\rVert_\infty$, and taking the continuous $L^2(0,L)$ norm gives
\begin{equation}
\lVert\widetilde u-u\rVert_{2}\le\frac{2\nu}{1-\lVert\varepsilon\rVert_\infty}\,\lVert\varepsilon_x\rVert_{2}.
\label{eq:colehopf-bound}
\end{equation}
Three consequences follow. A constant rescaling of $\phi$ cancels, which is why the normalization $C$ in~\eqref{eq:cole-hopf-initialization} is immaterial. The bound depends on the derivative of the relative error and on the denominator margin $1-\lVert\varepsilon\rVert_\infty$. For small relative perturbations, the first-order response is $-2\nu\varepsilon_x$. This is why the spatial structure of a perturbation matters in the amplification control of Section~\ref{sec:colehopf-amplification}. The white and smoothed perturbations there have equal absolute amplitude in the transformed field, and the measured ratio of about $12$ between their recovery errors is specific to that grid, field, and perturbation family. Finally, where $\phi$ is small the same absolute error in $\phi$ is a larger relative error $\varepsilon$, which is why the recovery is most sensitive near the shock center. The discrete recovery differs from~\eqref{eq:colehopf-perturbation} through the central difference, the kernel smoothing, and the positivity floor of Section~\ref{sec:burgers-benchmark}. The identity describes the continuous map that these approximate.

The Cole--Hopf transformation converts prescribed Burgers boundary values into derivative conditions on $\phi$~\citep{yan2023}. In the present unforced problem, setting $u(a,t)=u_b(t)$ in~\eqref{eq:cole-hopf-recovery} gives
\begin{equation}
\phi_x(a,t)=-\frac{u_b(t)}{2\nu}\,\phi(a,t),
\label{eq:colehopf-robin}
\end{equation}
a Robin condition on $\phi$ rather than a prescribed value of $\phi$. The fixed-endpoint treatment of Section~\ref{sec:burgers-benchmark} sets $\phi$ instead, which is why it carries the deterministic error of Table~\ref{tab:burgers-domain}. The exact time-dependent endpoint values used in the boundary control of Section~\ref{sec:colehopf-boundary} are available only because the transformed stationary shock is known in closed form, and they serve here as a diagnostic control. A particle treatment that imposes~\eqref{eq:colehopf-robin} directly is the natural next step.

\subsection{Scope of the diagnostics}
\label{sec:analysis-scope}
We state these four results for the one-dimensional calculations considered here. Their assumptions specify which parts of the diagnostic procedure transfer to other settings beyond the scope of this paper. The edge interpretation of a binned cumulative sum, the identity~\eqref{eq:cdf-identity}, the compatibility condition~\eqref{eq:boundary-compat}, and the recovery identity~\eqref{eq:colehopf-perturbation} each say what must be matched between a computation and its reference before a refinement curve is read. In two dimensions, gradient particles carry vector weights, as in the construction of Sherman and Mascagni~\citep{shermanmascagni1994}. Field recovery requires a potential or a line integral, while reflection depends on the domain geometry. A scalar Cole--Hopf potential gives $\boldsymbol u=-2\nu\nabla\log\phi$, which is irrotational. The diagnostic questions guide multidimensional studies, while the correction formulas require derivations appropriate to the reconstruction and geometry.

\section{Summary of Numerical Findings}
\label{sec:results-summary}

Table~\ref{tab:summary} collects the principal quantitative findings and their interpretations. The fitted rates are empirical and apply to the tested problems, parameter ranges, and implementations described above.

\begin{table}[pos=!htbp]
\centering
\small
\caption{Principal findings and interpretations. Heat and reaction--diffusion profile errors use $L_h^2$, while reaction--diffusion front location and speed errors retain their physical units. The Burgers domain, boundary, and decoupled-refinement values use RMSE.}
\label{tab:summary}
\begin{tabularx}{\textwidth}{@{}>{\raggedright\arraybackslash}p{0.26\textwidth} >{\raggedright\arraybackslash}p{0.37\textwidth} >{\raggedright\arraybackslash}X@{}}
\toprule
Study & Finding & Interpretation\\
\midrule
\addlinespace[0.4em]
\multicolumn{3}{@{}l}{\textit{Error attribution through controlled numerical comparisons}}\\
\addlinespace[0.4em]
Heat, coupled ensemble & spread $N^{-0.468}$, CI $[-0.521,-0.416]$ & sampling noise, compatible with $N^{-1/2}$\\[0.15cm]
\addlinespace[0.25em]
Heat, fixed bins, bin-center comparison & bias floors $4.6/3.6\times10^{-3}$ at $300/400$ bins & first-order evaluation error, predicted by the deterministic control to within 7\%\\[0.15cm]
\addlinespace[0.25em]
Heat, fixed bins, bin-edge comparison & bias $1.1\times10^{-3}$ at both bin counts & domain truncation error, and the apparent bin floor is removed\\[0.15cm]
\addlinespace[0.25em]
Heat, two-step prediction & predicted floors and transition counts of Table~\ref{tab:two-step-predictions}; measured/predicted $E_{\mathrm{total}}$ median $0.98$ over 160 comparisons & the floor and the particle count needed for a stated accuracy are computable in advance\\[0.15cm]
\addlinespace[0.25em]
Burgers, boundary control & $E_{\mathrm{det}}$, $0.172\to3.5\times10^{-4}$ with exact transformed data & fixed-endpoint treatment, reduced by exact boundary data\\[0.15cm]
\addlinespace[0.25em]
Burgers, decoupled refinement & particle RMSE $P^{-0.495}$, reconstruction-bin change about $1.3\%$, lowest tested particle RMSE at $\sigma_x=0.24$ & discrepancy from the fixed-endpoint reference decays as $P^{-0.495}$ under initialization-point refinement at fixed bins and bandwidth, interval containing $-1/2$\\[0.15cm]
\addlinespace[0.25em]
\midrule
\addlinespace[0.4em]
\multicolumn{3}{@{}l}{\textit{Verification and sensitivity}}\\
\addlinespace[0.4em]
Scalar reaction--diffusion ensemble & profile $N^{-0.460}$, front location and speed $N^{-0.524}$, aligned-profile $N^{-0.444}$ & shape and translation errors both decrease, location and profile rates compatible with $-1/2$, aligned-profile rate measurably shallower\\[0.15cm]
\addlinespace[0.25em]
Burgers, domain ensemble & $E_{\mathrm{det}}$ $0.172\!\to\!0.107$, $E_{\mathrm{GRW}}$ $0.104\!\to\!0.160$, total $0.187$--$0.200$, no monotone trend & composition shifts with domain, crossover between $L=6$ and $8$\\[0.15cm]
\addlinespace[0.25em]
Burgers, amplification control & white response about $12\times$ the smoothed at equal amplitude & recovery error depends on the spatial structure of the perturbation\\
\bottomrule
\end{tabularx}
\end{table}

\FloatBarrier

\section{Conclusions}
\label{sec:discussion}
Controlled comparisons identify reconstruction alignment and transformed boundary data as the sources of two apparent accuracy limits in gradient random walk calculations. For the heat equation, the cumulative bin sum represents the field at the bin edge. Comparing it there removes the first-order alignment error, and the remaining bias agrees with the domain truncation error. For Burgers' equation, exact transformed endpoint data remove the dominant deterministic error in the Cole--Hopf calculation. The particle runs retain their original fixed endpoints and approach the corresponding deterministic reference under initialization-point refinement.

The heat analysis also supports prediction before computation. For independent particles with fixed weights and known position distributions, the identity~\eqref{eq:weighted-identity} separates the expected squared error into a sampling term and a deterministic floor. Both terms are computable at the comparison points used in the reported norm. The two-step problem with unequal particle weights verifies the predicted floors and transitions across particle counts. A further ensemble at a count selected in advance agrees with its predicted error to within $1\%$. The count rule therefore provides a practical way to select a reconstruction convention and plan particle effort for an expected root-mean-square accuracy, with a margin appropriate to ensemble fluctuation.

For the scalar reaction--diffusion equation with the FitzHugh--Nagumo-type cubic, we verify the reaction-driven weight update combined with the diffusion step. Refinement improves the front's location, speed, and shape. Measuring the profile both before and after alignment separates translation from shape error, and shows that these quantities can have different empirical rates. Across the three problems, the sampling-dominated errors have fitted rates compatible with the Monte Carlo rate over the tested ranges, while deterministic floors explain the flattening of other refinement curves. The aligned reaction-front shape converges more slowly in the reported study. These distinctions are retained in the complete rate comparison in Table~\ref{tab:bootstrap}.

The analytical conditions in Section~\ref{sec:analysis} give each result its range of application. The cumulative-sum interpretation applies to signed weights, the expected-error prediction uses independent particles with specified fixed weights and known distributions, and the reaction update applies to differentiable scalar reaction laws with compatible boundary data. The Cole--Hopf identities describe both the transformed boundary relation and the sensitivity of recovery to the spatial variation of the relative error. These provide specific checks for other GRW calculations. The experiments demonstrate them in one-dimensional heat, reaction--diffusion, and Burgers problems.

A natural extension is to impose the transformed boundary relation directly in the Burgers particle solver. Further studies can examine smoothing choices across shock amplitudes and viscosities, and adapt the comparisons to multidimensional reconstruction. For the present calculations, the practical conclusion is direct. Match the reconstructed quantity and its comparison points, impose boundary data consistent with the physical problem, and use the resulting error prediction to guide refinement wherever the particle distribution is known.\clearpage
\appendix
\renewcommand{\thesection}{\Alph{section}}
\numberwithin{equation}{section}
\section{\texorpdfstring{Identities for the canonical problems}{Identities for the canonical problems}}
\label{app:benchmarks}

For the logistic front~\eqref{eq:fhn-exact}, direct differentiation gives
\begin{equation}
u_x=\tfrac12u(1-u),\qquad
u_{xx}=\tfrac14u(1-u)(1-2u),\qquad
u_t=\tfrac\theta2u(1-u).
\end{equation}
Substituting into $u_t-Du_{xx}=f(u)$ yields~\eqref{eq:fhn-reaction}. Differentiating that cubic yields~\eqref{eq:fhn-statistic}. At $t=0$, solving $u^{\mathrm{ex}}(X_i^0,0)=q_i$ gives the quantile initialization~\eqref{eq:fhn-quantile-initialization}.

For the Burgers shock, set $z=A(x-x_c)/(2\nu)$. Then
\begin{equation}
u=-A\tanh z,\quad u_x=-\frac{A^2}{2\nu}\operatorname{sech}^2z,\quad
u_{xx}=\frac{A^3}{2\nu^2}\operatorname{sech}^2z\tanh z,
\end{equation}
so that
\begin{equation}
u u_x=\nu u_{xx},\qquad u_t=0,
\label{eq:app-shock-stationary}
\end{equation}
and the profile is a stationary solution of~\eqref{eq:burgers-pde}. The far-field states $u_L=A$, $u_R=-A$ also give the Rankine--Hugoniot speed $(u_L+u_R)/2=0$ for Burgers' flux~\citep{leveque2002}.

With $k=A/(2\nu)$, the transformed initialization follows from
\begin{equation}
\Psi_0(x)=-2\nu\log\frac{\cosh(k(x-x_c))}{\cosh(kx_c)}.
\end{equation}
For the centered shock $x_c=L/2$, normalization to $\max_x\phi_0=1$ gives
\begin{equation}
\phi_0(x)=\frac{\cosh\big(k(x-x_c)\big)}{\cosh(kx_c)}.
\label{eq:app-phi0}
\end{equation}
Since $\phi_0''=k^2\phi_0$, the heat solution is
\begin{equation}
\phi(x,t)=e^{\nu k^2t}\phi_0(x).
\label{eq:app-phi-growth}
\end{equation}
Its growing endpoint values and stationary recovered Burgers field are therefore compatible.

\section{Heat uncertainty and reconstruction details}
\label{app:heat-uncertainty}

\subsection{\texorpdfstring{Reflected particle distribution and empirical covariance}{Reflected particle distribution and empirical covariance}}
Let $\gamma_t(z)=(4\pi\alpha t)^{-1/2}\exp[-z^2/(4\alpha t)]$, with $\alpha>0$, $t>0$, and $x_0\in(0,L)$. Folding a Gaussian position on the real line with mean $x_0$ and variance $2\alpha t$ into $[0,L]$ gives the density
\begin{equation}
p_R(x,t\mid x_0)=\sum_{n\in\mathbb Z}\big[\gamma_t(x-x_0-2nL)+\gamma_t(x+x_0-2nL)\big],\qquad 0<x<L.
\label{eq:reflected-density}
\end{equation}
The image pairs are even about each wall, so the spatial derivative vanishes at $0$ and $L$. The inverse images of $(0,L)$ partition the real line up to endpoints, giving unit mass, and each Gaussian solves the heat equation. Thus this is the classical Neumann heat kernel~\citep{carslawjaeger1959}, equivalently the transition density of reflected Brownian motion~\citep{pavliotis2014}. The distribution function of a reflected particle at the final time is therefore
\begin{equation}
F_R(x)=\int_0^x p_R(y,T\mid x_0)\,dy .
\label{eq:reflected-cdf}
\end{equation}
The Gaussian antiderivative is $\frac12[1+\operatorname{erf}(z/\sqrt{4\alpha T})]$, which yields the difference $G(x)-G(0)$ used in Section~\ref{sec:heat-grid-paired}. The finite image range of $G$ is a numerical approximation to the infinite sum, with a neglected tail below machine precision for the reported heat parameters.

The implemented walk folds the position into $[0,L]$ after every step, and its distribution equals that of the unreflected Gaussian walk folded once at the end. Let $\Pi$ be the folding map, which is even and $2L$-periodic, so that $\Pi(y)=\epsilon y+2nL$ for some sign $\epsilon$ and integer $n$ depending on $y$. For any real $y$ and any increment $\xi$ whose distribution is symmetric,
\begin{equation}
\Pi\big(\Pi(y)+\xi\big)=\Pi(\epsilon y+2nL+\xi)=\Pi(\epsilon y+\xi)=\Pi(y+\epsilon\xi),
\label{eq:folding-identity}
\end{equation}
and $\epsilon\xi$ has the same distribution as $\xi$. Hence $\Pi(\Pi(y)+\xi)$ and $\Pi(y+\xi)$ have the same distribution. Now let
\[
Y_n=x_0+\sum_{k\le n}\xi_k,\qquad
X_n=\Pi(X_{n-1}+\xi_n),\qquad
X_0=Y_0=x_0,
\]
be the unreflected Gaussian walk on the real line and the implemented reflected walk, respectively, with independent increments $\xi_k\sim\mathcal N(0,2\alpha\Delta t)$. If $X_{n-1}$ has the distribution of $\Pi(Y_{n-1})$, the identity gives $X_n$ the distribution of $\Pi(Y_{n-1}+\xi_n)=\Pi(Y_n)$, so by induction $X_K$ has the distribution of $\Pi(Y_K)$. Since $Y_K$ is Gaussian with mean $x_0$ and variance $2\alpha T$, folding it gives~\eqref{eq:reflected-density}.

For the indicator variables in Section~\ref{sec:analysis-cdf}, $I_i(x)I_i(y)=I_i(\min(x,y))$. Independence across particles therefore gives the exact covariance
\begin{equation}
\operatorname{Cov}(F_N(b_j),F_N(b_k))=\frac{\Gamma_{jk}}{N},\qquad
\Gamma_{jk}=F(\min(b_j,b_k))-F(b_j)F(b_k).
\label{eq:cdf-covariance}
\end{equation}
Here $F=F_R$. Correlations between evaluation points matter for uncertainty in an error norm even though the expected squared norm in~\eqref{eq:cdf-identity} uses only the diagonal entries.

\subsection{Exact finite-ensemble moments}
For $S\ge2$ independent, identically distributed realizations with finite second moments, let $b=\mathbb E[u_{N,s}]-u^{\mathrm{ref}}$ and set
\[
Z_s=u_{N,s}-\mathbb E[u_{N,s}],\qquad
\bar Z=\frac1S\sum_{s=1}^{S} Z_s,\qquad
V=\mathbb E\lVert Z_s\rVert_{L_h^2}^2 .
\]
Independence gives $\mathbb E\lVert\bar Z\rVert_{L_h^2}^2=V/S$. The identity
\[
\sum_s\lVert Z_s-\bar Z\rVert_{L_h^2}^2=\sum_s\lVert Z_s\rVert_{L_h^2}^2-S\lVert\bar Z\rVert_{L_h^2}^2
\]
then yields the first statement below, and expanding $\lVert b+\bar Z\rVert_{L_h^2}^2$ with $\mathbb E\bar Z=0$ gives the second,
\begin{equation}
\mathbb E[E_{\mathrm{spread}}^2]=\frac{S-1}{S}V,\qquad
\mathbb E[E_{\mathrm{bias}}^2]=\lVert b\rVert_{L_h^2}^2+\frac VS ,
\label{eq:finite-ensemble-factors}
\end{equation}
as used in Section~\ref{sec:ensemble-errors}.

For heat, define $\Sigma=h\Gamma/N$, $q_j=\sqrt h\,[F(b_j)-g_j]$, and $\bar\Sigma=\Sigma/S$. The centered reconstruction vector multiplied by $\sqrt h$ has covariance $\Sigma$. Its ensemble mean has covariance $\bar\Sigma$. Hence
\begin{equation}
\mathbb E[E_{\mathrm{bias}}^2]=q^Tq+\operatorname{tr}\bar\Sigma,\qquad
\mathbb E[E_{\mathrm{spread}}^2]=\frac{S-1}{S}\operatorname{tr}\Sigma.
\label{eq:heat-exact-moments}
\end{equation}
These are exact second-moment statements for the independent-particle model. We compare the measured spreads with the root-mean-square value $\sqrt{\mathbb E[E_{\mathrm{spread}}^2]}$. At $N=2000$ and $5000$ the measured spreads agree with it to within $0.4\%$, at $N=500$ and $1000$ to within $4\%$, and at $N\ge10{,}000$ they lie $11$--$19\%$ above it.

\subsection{Gaussian and delta approximations}
Approximate the centered empirical fields by independent Gaussian vectors with covariance $\Sigma$, so the ensemble mean $z$ has covariance $\bar\Sigma=\Sigma/S$. This distributional assumption is additional to the exact moments in~\eqref{eq:heat-exact-moments}. The Gaussian quadratic-form variance formula~\citep{rencher_schaalje2007quadratic}, applied to $\lVert q+z\rVert^2$, gives
\begin{equation}
\operatorname{Var}(E_{\mathrm{bias}}^2)\approx4q^T\bar\Sigma q+2\operatorname{tr}(\bar\Sigma^2).
\label{eq:heat-bias-variance-approx}
\end{equation}
Under the same model,
\[
E_{\mathrm{spread}}^2\ \sim\ \frac1S\sum_j\lambda_j\,\chi^2_{S-1,j},
\]
with independent chi-squared variables and eigenvalues $\lambda_j$ of $\Sigma$. Since $\operatorname{Var}(\chi^2_{S-1})=2(S-1)$,
\begin{equation}
\operatorname{Var}(E_{\mathrm{spread}}^2)\approx\frac{2(S-1)}{S^2}\operatorname{tr}(\Sigma^2).
\label{eq:heat-spread-variance-approx}
\end{equation}
For either squared statistic $Q$ with $\mathbb EQ>0$, the first-order delta method for $\sqrt Q$~\citep{vandervaart1998} gives
\begin{equation}
\operatorname{sd}(\sqrt Q)\approx\frac{\sqrt{\operatorname{Var}(Q)}}{2\sqrt{\mathbb EQ}}.
\label{eq:heat-delta-approx}
\end{equation}
This requires fluctuations small enough for linearization and estimates the relative standard deviation of the spread as approximately $7\%$. These approximations give the scale of the fluctuations of the measured norms.

To check the excess at the larger counts, we repeated the $N=50{,}000$, $M=300$ computation with three further blocks of thirty seeds, starting at 1000, 2000, and 3000. The original seeds 42--71 reproduce the spread of $3.30\times10^{-3}$ in Table~\ref{tab:heat-grid-paired}, and the three further blocks give spreads between $2.68$ and $2.81\times10^{-3}$, which bracket the root-mean-square value $2.77\times10^{-3}$. This supports reading the excess as a fluctuation of the original block of seeds.

One detail of the binning is worth recording. The implementation assigns a particle lying exactly on a bin edge to the bin on its right, whereas~\eqref{eq:cumulative-bins} counts it in the bin on its left. Because the particle positions have a continuous distribution at any positive time, a particle lies exactly on an edge with probability zero, so the two conventions give the same reconstruction for the heat computations. The convention matters only for deterministic test inputs placed exactly on edges.

\subsection{\texorpdfstring{Weighted particles and the particle-count rule}{Weighted particles and the particle-count rule}}
\label{app:weighted}
Let two groups of globs be released at $a$ and $b$ with deterministic weights $A/N_a$ and $B/N_b$ and independent increments, and let $I_{a,i}(x)=\mathbf 1_{\{X_{a,i}\le x\}}$ and $I_{b,i}(x)=\mathbf 1_{\{X_{b,i}\le x\}}$. The bin-and-sum reconstruction at a right edge $b_j$ is
\[
U_N(b_j)=u_L+\frac{A}{N_a}\sum_{i=1}^{N_a}I_{a,i}(b_j)+\frac{B}{N_b}\sum_{i=1}^{N_b}I_{b,i}(b_j).
\]
Each indicator is a Bernoulli variable with success probability $F_a(b_j)$ or $F_b(b_j)$, and the indicators of distinct globs are independent, so at each edge
\begin{equation}
\mathbb E[U_N(b_j)]=m(b_j),\qquad
\operatorname{Var}[U_N(b_j)]=\frac{A^2}{N_a}F_a(1-F_a)+\frac{B^2}{N_b}F_b(1-F_b),
\label{eq:weighted-mean-variance}
\end{equation}
with the distribution functions evaluated at $b_j$. For the covariance between two edges, apply the indicator argument of~\eqref{eq:cdf-covariance} separately to $F_a$ and $F_b$, multiply the two contributions by $A^2/N_a$ and $B^2/N_b$, and add them. The norm-weighted covariance includes a further factor $h$. Omitting the second group and setting $N_a=N$ and $A=u_R-u_L$ recovers the single-step case of Section~\ref{sec:analysis-cdf}. Applying $\mathbb E[(Y-g)^2]=\operatorname{Var}(Y)+(\mathbb EY-g)^2$ at each edge, summing with weight $h$, and substituting $N_a=f_aN$ and $N_b=f_bN$ gives~\eqref{eq:weighted-identity}. The finite-ensemble factors then follow from~\eqref{eq:finite-ensemble-factors} with $V=V_h/N$,
\begin{equation}
\mathbb E[E_{\mathrm{total}}^2]=B_h^2+\frac{V_h}{N},\qquad
\mathbb E[E_{\mathrm{bias}}^2]=B_h^2+\frac{V_h}{NS},\qquad
\mathbb E[E_{\mathrm{spread}}^2]=\frac{S-1}{S}\cdot\frac{V_h}{N}.
\label{eq:weighted-ensemble-factors}
\end{equation}

These statements were checked against the implementation. Drawing four thousand ensembles of twelve realizations directly from the final-time distribution at $M=200$ and $N=400$ reproduces the three moments of~\eqref{eq:weighted-ensemble-factors} to within one standard error, and the sampling coefficient obtained from the exact covariance agrees with the value in Table~\ref{tab:two-step-predictions}.

Equating the sampling and deterministic terms gives the transition count~\eqref{eq:transition-count}. Solving $B_h^2+V_h/N\le\tau^2$ for $N$, when $\tau>B_h$, gives the particle-count rule~\eqref{eq:particle-count-rule} used in Section~\ref{sec:heat-two-step}.

\subsection{\texorpdfstring{Resampling scheme for the fitted trends}{Resampling scheme for the fitted trends}}
\label{app:bootstrap}
Each trend study in Table~\ref{tab:bootstrap} reuses one seed list across its parameter values, so realizations at different values reuse pseudorandom number sequences. Table~\ref{tab:bootstrap} reports each fitted trend under the two resampling schemes of Section~\ref{sec:ensemble-bootstrap}: resampling the seed labels jointly across parameter values, which preserves that grouping and is the scheme quoted in the text, and resampling realizations independently within each value. The joint intervals are the wider of the two in all but two cases, the front-location and initialization-point trends. The two schemes agree on which trends are compatible with $-1/2$ except for the heat total errors of the coupled and bin-edge treatments and the reaction profile error, which the independent scheme places just above it. The heat entries use the realization profiles of the paired study, which reproduce the statistics of Table~\ref{tab:heat-grid-paired} to the stored precision.

\begin{table}[pos=htbp]
\centering\small
\caption{Fitted trends under joint and independent resampling. Log--log slopes except the Burgers domain trend, which is a linear fit against $L$.}
\label{tab:bootstrap}
\begin{tabular}{@{}llccc@{}}
\toprule
Study & Statistic & Fit & Joint by seed & Independent\\
\midrule
Heat, coupled & total & $-0.457$ & $[-0.505,-0.410]$ & $[-0.493,-0.423]$\\
Heat, coupled & spread & $-0.468$ & $[-0.521,-0.416]$ & $[-0.508,-0.430]$\\
Heat, coupled & bias & $-0.262$ & $[-0.461,-0.201]$ & $[-0.451,-0.209]$\\
Heat, fixed 300, bin center & total & $-0.353$ & $[-0.395,-0.313]$ & $[-0.388,-0.319]$\\
Heat, fixed 300, bin center & spread & $-0.468$ & $[-0.522,-0.416]$ & $[-0.508,-0.430]$\\
Heat, fixed 300, bin center & bias & $-0.053$ & $[-0.225,+0.075]$ & $[-0.200,+0.056]$\\
Heat, fixed 300, bin edge & total & $-0.456$ & $[-0.503,-0.410]$ & $[-0.492,-0.422]$\\
Heat, fixed 300, bin edge & bias & $-0.229$ & $[-0.447,-0.184]$ & $[-0.433,-0.202]$\\
Reaction front & profile & $-0.460$ & $[-0.505,-0.415]$ & $[-0.497,-0.423]$\\
Reaction front & location and speed & $-0.524$ & $[-0.605,-0.445]$ & $[-0.609,-0.442]$\\
Reaction front & aligned profile & $-0.444$ & $[-0.479,-0.409]$ & $[-0.474,-0.412]$\\
Burgers, initialization points & particle RMSE & $-0.495$ & $[-0.521,-0.470]$ & $[-0.523,-0.466]$\\
Burgers, domain size & total RMSE vs $L$ & $-1.20\times10^{-3}$ & $[-3.02,+0.73]\times10^{-3}$ & $[-2.90,+0.62]\times10^{-3}$\\
\bottomrule
\end{tabular}
\end{table}
\FloatBarrier

\Needspace{12\baselineskip}
\section*{Code and Data Availability Statement}

The code, input data, and reference results needed to reproduce the numerical studies and figures are archived on Zenodo as version~1.1.0 at \url{https://doi.org/10.5281/zenodo.22680376}~\citep{abkin_daripa_grw_2026}.

\section*{Acknowledgments}

The authors thank Oliver Stalker for providing an early version of the Python code.

\bibliographystyle{plainnat} % numbered in alphabetical order of authors
\bibliography{references}

\end{document}